\documentclass[reqno]{amsart}
\usepackage{graphicx,amsmath,amssymb,hyperref,nicefrac,microtype,xcolor,bm,enumitem}
\usepackage[foot]{amsaddr}
\usepackage[paperwidth=7in,paperheight=10in,text={5.4in,7.95in},centering,top=1in]{geometry}

\usepackage[numbers,square,comma]{natbib}

\usepackage{amsthm}
\theoremstyle{plain}
\newtheorem{theorem}{Theorem}
\newtheorem{algorithm}[theorem]{Algorithm}

\numberwithin{equation}{section}

\numberwithin{theorem}{section}

\usepackage[T1]{fontenc}
\usepackage{mathptmx}

\newcommand{\ignore}[1]{}
\newcommand{\opn}{\operatorname}
\newcommand{\veps}{\varepsilon}
\newcommand{\re}{\operatorname{Re}}
\newcommand{\im}{\operatorname{Im}}

\newcommand{\mbb}[1]{\mathbb{#1}}

\newcommand{\mc}[1]{\mathcal{#1}}
\newcommand{\pa}{\partial}
\newcommand{\der}[2]{\frac{\partial #1}{\partial #2}}

\newcommand{\jd}{\displaystyle}

\newcommand{\jt}{\textstyle}

\newcommand{\wdg}[1]{{#1}^{\scriptscriptstyle\bm\wedge}}
\newcommand{\e}[1]{{(#1)}}

\newcommand{\gridkap}{\kappa}

\newcommand{\la}{\langle}
\newcommand{\ra}{\rangle}
\newcommand{\lkr}{\langle k\rangle}
\newcommand{\hh}{h}

\begin{document}

\title{Harmonic Stability of Standing Water Waves}

\author{Jon Wilkening}
\email{wilken@math.berkeley.edu}
\thanks{Dedicated to Walter Strauss in honor of his 80th birthday.
  This work was supported in part by the National Science Foundation under
  award number DMS-1716560 and by the Department of Energy, Office of
  Science, Applied Scientific Computing Research, under award number
  DE-AC02-05CH11231.}
\address{Department of Mathematics, University of
  California, Berkeley, CA 94720-3840}


\begin{abstract}
  A numerical method is developed to study the stability of standing
  water waves and other time-periodic solutions of the free-surface
  Euler equations using Floquet theory.  A Fourier truncation of the
  monodromy operator is computed by solving the linearized Euler
  equations about the standing wave with initial conditions ranging
  over all Fourier modes up to a given wave number. The eigenvalues of
  the truncated monodromy operator are computed and ordered by the
  mean wave number of the corresponding eigenfunctions, which we
  introduce as a method of retaining only accurately computed Floquet
  multipliers. The mean wave number matches up with analytical results
  for the zero-amplitude standing wave and is helpful in identifying
  which Floquet multipliers collide and leave the unit circle to form
  unstable eigenmodes or rejoin the unit circle to regain stability.
  For standing waves in deep water, most waves with crest acceleration
  below $A_c=0.889$ are found to be linearly stable to harmonic
  perturbations; however, we find several bubbles of instability at
  lower values of $A_c$ that have not been reported previously in the
  literature. We also study the stability of several new or recently
  discovered time-periodic gravity-capillary or gravity waves in deep
  or shallow water, finding several examples of large-amplitude waves
  that are stable to harmonic perturbations and others that are not.
  A new method of matching the Floquet multipliers of two nearby
  standing waves by solving a linear assignment problem is also
  proposed to track individual eigenvalues via homotopy from the
  zero-amplitude state to large-amplitude standing waves.
\end{abstract}


\maketitle

\markboth{JON WILKENING}{HARMONIC STABILITY OF STANDING WATER WAVES}

\section{Introduction}

Standing water waves have played a central role in the study
of fluid mechanics throughout history. Already in 1831, Faraday
observed patterns of ink at the surface of milk driven by a tuning
fork \cite{faraday}. He noticed that the waves oscillate at half the
forcing frequency, which was later explained by Rayleigh
\cite{rayleigh:faraday}. Benjamin and Ursell performed a linear
stability analysis in \cite{benjamin:ursell}, which was later extended
to viscous fluids by Kumar and Tuckerman \cite{kumar94}. Additional
theoretical, numerical and experimental studies of Faraday waves
include \cite{vinyals97, skeldon07, knobloch:faraday,
  wright00, chen02, murakami01, oconnor08, perinet, henderson91,
jiang:ting, bush}.

Standing waves at the surface of an inviscid fluid with no external
driving force have also been studied extensively. Building on the
asymptotic techniques developed by Stokes for traveling waves
\cite{stokes:1880}, Rayleigh showed how to incorporate time-evolution
in the analysis and computed standing waves to third order
\cite{rayleigh:1915}. Penney and Price extended Rayleigh's expansion
to 5th order and conjectured the existence of an extreme standing wave
that forms a 90 degree corner each time the wave comes to rest
\cite{penney:52}. Taylor confirmed experimentally that large-amplitude
standing waves nearly form 90 degree corners \cite{taylor:53} but was
skeptical of Penney and Price's mathematical argument.  Tadjbakhsh and
Keller \cite{tadjbakhsh} and Concus \cite{concus:62} incorporated
finite-depth and surface tension effects in the asymptotic expansions.
Concus later realized there were small-divisor issues that called into
question the validity of the asymptotic techniques that had been
developed to study standing waves up to that point
\cite{concus:64}. It took 40 years for these issues to be sorted out
using Nash-Moser theory \cite{plotnikov01, iooss05, alazard:14}. Along
the way, asymptotic expansions to arbitrary order were computed for
standing waves in deep water by Schwartz and Whitney
\cite{schwartz:81} and a semi-analytic theory of standing waves was
proposed by Amick and Toland \cite{amick:87}.  On the computational
side, Mercer and Roberts \cite{mercer:92, mercer:94}, Tsai and Jeng
\cite{tsai:jeng:94}, Bryant and Stiassnie \cite{bryant:stiassnie:94},
Schultz et al.~\cite{schultz}, Okamura \cite{okamura:03, okamura:10},
Wilkening \cite{water1}, Wilkening and Yu \cite{water2}, and many
others have developed numerical algorithms for computing standing
water waves. Due to the difficulty of maintaining accuracy in nearly
singular free-surface flow calculations, different conclusions were
reached by these authors regarding the form of the largest-amplitude
standing wave. The question was finally settled in
\cite{water1,water2}, where it was shown that the Penney-Price
conjecture is false due to a breakdown in self-similarity at the
crests of very large-amplitude standing waves.

While much is known about the stability of traveling water waves
\cite{benjamin67, benjamin67b, lh:stab:78a, lh:stab:78b, crawford:81,
  mackay:86, ioualalen:93, nicholls09, oliveras11, wilton:16},
relatively little work has been done on the stability of standing
waves. The most comprehensive study to date is by Mercer and Roberts
\cite{mercer:92}, who found that standing waves on deep water are
linearly stable to harmonic perturbations for crest acceleration in
the range $0\le A_c\le0.889$.  (Crest acceleration is the downward
  acceleration of a fluid particle at the wave crest at maximum height
  when the fluid comes to rest, relative to the gravitational
  acceleration.)  They also studied subharmonic perturbations with
wavelength equal to 8 times the fundamental wavelength of the standing
wave, finding sideband instabilities in which the base wave (mode
  $m=8$) interacts with modes $m\pm l$, with $l$ an integer.  Bryant
and Stiassnie \cite{bryant:stiassnie:94} found that sideband
instabilities lead to cyclic recurrence over hundreds of periods of
the basic standing wave when evolved numerically on a domain
containing 9 replicas of the standing wave.  Bridges and Laine-Pearson
\cite{bridges:counter:prop:2004} have shown that the modulational
instability of standing waves is closely linked to that of the
component counter-propagating traveling waves at the weakly nonlinear
level for a wide range of Hamiltonian PDEs, including water waves.
In the present work we focus on harmonic stability, which is natural
in the context of standing waves in a container with vertical walls.
Subharmonic stability of more general ``traveling-standing'' waves will
be considered elsewhere \cite{waterTS}.

To determine the stability of standing waves and other time-periodic
solutions of the free-surface Euler equations to harmonic
perturbations, we compute Floquet multipliers of the monodromy
operator to determine if they lie on the unit circle.  The main
difference between our approach and that of Mercer and Roberts
\cite{mercer:92} is that we employ a Fourier basis rather than
discrete delta functions.  This has a computational advantage in that
only the leading columns of the monodromy operator need be computed;
thus, the mesh can be refined independently of the size of the
(truncated) operator so that all the matrix entries are accurate.  We
also order the eigenvalues according to the ``mean wave number'' of
the corresponding eigenfunctions, which turns out to be closely
correlated with the residual error of the eigenpair.  By parallelizing
the computation of the operator using a GPU, we are able to increase
the number of Floquet multipliers that can be computed by 2 orders of
magnitude over previous studies.  Doing this reveals additional
bubbles of instability \cite{mackay:86} in the range $0\le A_c\le
0.889$, previously thought to contain only stable solutions.

In addition to studying the stability of standing waves in deep water,
we consider for the first time finite depth and capillary effects on
stability.  In infinite depth, standing waves involve large-scale
motion of the bulk fluid. However, in shallow water, they are better
described as counter-propagating solitary waves that repeatedly
collide with one another.  We consider a particular family of standing
waves with wavelength $2\pi$ and fluid depth $h=0.05$. We find
examples of waves well outside of the weakly nonlinear regime that are
stable to harmonic perturbations, and another that is unstable. We
also study the stability of gravity-capillary waves. These can take
the form of counter-propagating depression waves, for which we give an
example of a large-amplitude wave that is stable to harmonic
perturbations.  We also search for a new type of time-periodic water
wave consisting of counter-propagating gravity-capillary solitary
waves \cite{lh89,vandenBroeck:cap, milewski:11} that collide
repeatedly with each other.  We present two such waves, one that is
stable to harmonic perturbations and one that is unstable.

A technical challenge that arises in studying the stability of
standing waves is that the Floquet multipliers are complex numbers
whose phase is only known modulo $2\pi$; thus, while they mostly
lie on the unit circle, they are scrambled together. By contrast, for traveling waves,
the eigenvalues mostly lie on the imaginary axis and are naturally
ordered by their imaginary part. We introduce a ``mean wave number''
to measure how oscillatory the corresponding eigenfunctions are, and
use this to order the eigenvalues. This aids in discarding inaccurate
eigenvalues and also proves useful in identifying which eigenvalues
collide and leave the unit circle to form unstable eigenmodes or
rejoin the unit circle to regain stability. In order to track
individual eigenvalues via a homotopy method, we propose an algorithm
to match the Floquet multipliers of nearby standing waves by solving a
linear assignment problem. The mean wave number also plays a role in
the cost matrix of this assignment problem.

This paper is organized as follows. In Section~\ref{sec:osm} we
discuss the free-surface Euler equations, non-uniform meshes,
symmetry, and our overdetermined shooting method for computing
time-periodic solutions. In Section~\ref{sec:stable}, we discuss
linearization about time-periodic solutions, Hamiltonian time
reversal, and a Fourier basis for computing the monodromy operator
efficiently. We also analyze the stability of the zero-amplitude
solution and present our numerical algorithm for computing Floquet
multipliers. In Section~\ref{sec:deep} we study the stability of
standing waves in deep water with crest acceleration in the range
$0\le A_c\le 0.9639$.  Like Mercer and Roberts \cite{mercer:92}, we
find that there is a critical bifurcation at $A_c=0.889$ beyond which
all standing waves are unstable. However, we also find new bubbles of
instability below this critical threshold and investigate one of them
in detail. In Section~\ref{sec:jchains} we discuss the numerical
splitting of the $\lambda=1$ eigenvalue due to Jordan chains
associated with translation in time and space. In
Section~\ref{sec:shallow}, we study the stability of
counter-propagating solitary waves in shallow water, including
large-amplitude time-periodic waves that temporarily form a jet that
is taller than the fluid depth. In Section~\ref{sec:surf}, we study
the stability of two types of gravity-capillary waves in deep
water. The first type consists of counter-propagating depression waves
that bifurcate from the flat rest state as predicted by
Concus~\cite{concus:62}. The second is a new type of time-periodic
wave constructed from two counter-propagating solitary waves
\cite{lh89,vandenBroeck:cap,milewski:11} with initial conditions tuned
to achieve time-periodicity. The appendices describe (\ref{sec:BI})
the boundary integral formulation; (\ref{sec:trav}) a method of
computing traveling waves using the overdetermined shooting method;
(\ref{sec:J}) computation of the Jacobian and the state transition
matrix; and (\ref{sec:matching}) a linear assignment problem for
matching eigenvalues at adjacent values of $A_c$ to track eigenvalues
via homotopy from the zero-amplitude state to large-amplitude standing
waves.

\section{Preliminaries: computation of standing water waves}
\label{sec:osm}

In this section, we describe an overdetermined shooting algorithm for
computing time-periodic solutions of the free-surface Euler equations
with spectral accuracy. We build on this method in
Section~\ref{sec:stable} to compute the monodromy operator and its
eigenvalues, which determine whether the underlying solution is linearly
stable to harmonic perturbations. We also use the method in
Section~\ref{sec:surf} to compute new families of time-periodic
gravity-capillary waves and investigate their stability.  The
overdetermined shooting algorithm is explained in more detail in
Wilkening and Yu \cite{water2}, and builds on previous shooting
methods \cite{mercer:92, mercer:94, schultz, smith:roberts:99}.  Other
successful approaches for computing standing waves include Fourier
collocation in space and time \cite{vandenBroeck:81, tsai:jeng:94,
  okamura:03, ioualalen:03, okamura:10} and semi-analytic series
expansions \cite{schwartz:81,amick:87}.

\subsection{Equations of motion}
\label{sec:eqm}

The equations of motion of a free surface $\eta(x,t)$ evolving over
an ideal fluid with velocity potential $\phi(x,y,t)$ may be
written \cite{whitham74,johnson97,craik04,craik05}
\begin{equation}  \label{eq:ww}
\begin{aligned}
    \eta_t &= \phi_y - \eta_x\phi_x, \\[-3pt]
\varphi_t &= P\left[\phi_y\eta_t - \frac{1}{2}\phi_x^2 -
  \frac{1}{2}\phi_y^2 - g\eta +
  \frac{\sigma}{\rho}\partial_x\left(\frac{\eta_x}{\sqrt{1+\eta_x^2}}
  \right)\right],
\end{aligned}
\end{equation}
where subscripts denote partial derivatives, $\varphi(x,t) =
\phi(x,\eta(x,t), t)$ is the restriction of $\phi$ to the free
surface, $g$ is the acceleration of gravity, $\rho$ is the fluid
density, $\sigma\ge0$ is the surface tension (possibly zero), and $P$
is the $L^2$-projection to zero mean that annihilates constant
functions,
\begin{equation}
  P = \opn{id} - P_0, \qquad
  P_0 f = \frac{1}{2\pi}\int_0^{2\pi} f(x)\,dx.
\end{equation}
Here we have non-dimensionalized the equations so that the wavelength
of the standing wave is $2\pi$. In physical variables, (\ref{eq:ww})
still describes the dynamics if we put tildes over $\eta$, $\phi$,
$\varphi$, $x$, $y$, $t$, $g$, $\sigma$ and $\rho$ and modify $P_0$ to
compute the mean over $[0,2\pi L]$, where
\begin{equation}\label{eq:nondim}
  \tilde \eta(\tilde x,\tilde t) = L \eta\left(\frac{\tilde x}{L},
    \frac{\tilde t}{\tau}\right), \quad
  \tilde \phi(\tilde x,\tilde y,\tilde t) =
  \frac{L^2}{\tau} \phi\left(\frac{\tilde x}{L},\frac{\tilde y}{L},
    \frac{\tilde t}{\tau}\right), \quad
  \tilde g = \frac {L}{\tau^2}g, \quad
  \frac{\tilde\sigma}{\tilde\rho} = \frac{L^3}{\tau^2}\frac{\sigma}\rho.
\end{equation}
Here $L$ and $\tau$ are the physical length and time corresponding
to one unit of dimensionless length and time, respectively. 
The projection $P$ yields a convenient convention for selecting the
arbitrary additive constant in the potential, and has no effect
(barring roundoff errors and discretization errors) if the fluid has
infinite depth and the mean surface height is zero. The velocity
components $u=\phi_x$ and $v=\phi_y$ at the free surface can be
computed from $\varphi$ via
\begin{equation}\label{eq:uv:from:G}
  \begin{pmatrix} \phi_x \\ \phi_y \end{pmatrix} =
  \frac{1}{1+\eta'(x)^2}\begin{pmatrix}
    1 & -\eta'(x) \\ \eta'(x) & 1 \end{pmatrix}
  \begin{pmatrix}
    \varphi'(x) \\
    \mc G\varphi(x)
  \end{pmatrix},
\end{equation}
where a prime denotes a derivative and $\mc G$ is the
Dirichlet-Neumann operator \cite{craig:sulem:93}
\begin{equation}\label{eq:DNO:def}
  \mc G\varphi(x) = \sqrt{1+\eta'(x)^2}\,\, \der{\phi}{n}(x+i\eta(x))
\end{equation}
for the Laplace equation, with periodic boundary conditions in $x$,
$\phi=\varphi$ on the upper boundary, and $\phi_y=0$ at $y=-\infty$ or
at $y=-h$ in the finite depth case. We have suppressed $t$ in the
notation as time is frozen when solving the Laplace equation.  We
compute $\mc G\varphi$ using a boundary integral collocation method;
see Appendix~\ref{sec:BI}.

This method is flexible enough to allow for non-uniform meshes in
space and time. We adopt a simple approach in which time is divided
into $\nu$ segments of relative size $\theta_l$, where
$\theta_1+\theta_2+\cdots+\theta_\nu=1$. We take $N_l$ uniformly
spaced timesteps on segment $l$. The spatial grid on this time segment
is fixed, with two parameters $\gridkap_l\in\mathbb{Z}$ and
$\rho_l\in(0,1]$ controlling the relative grid spacing near $x=0$ and
$x=\pi$.  The spatial grid on segment $l$ is
$x_i=\xi_l(\alpha_i)$, where $\alpha_i=2\pi i/M_l$ for $0\le i<M_l$
and
\begin{equation}\label{eq:xi:def}
  \xi_l(\alpha) = \int_0^\alpha E_l(\beta)\,d\beta, \qquad\quad
  E_l(\alpha) = \left\{\begin{array}{rl}
  \!\!1-A_l\big[\sin^{2\gridkap_l}(\alpha/2)-\mu_l\big], & \gridkap_l>0, \!\!\\[2pt]
  \!\!1, & \gridkap_l=0,\!\!\\
  \!\!1-A_l\big[\cos^{2|\gridkap_l|}(\alpha/2)-\mu_l\big], & \gridkap_l<0.  \!\!
  \end{array}\right.
\end{equation}
Here\, $\jd\mu_l = \frac{(2|\gridkap_l|-1)!!}{(2|\gridkap_l|)!!} =
\frac{\Gamma(|\gridkap_l|+\frac12)}{\sqrt\pi\,\,\Gamma(|\gridkap_l|+1)}$
and\, $\jd A_l = \frac{1-\rho_l}{1-\mu_l(1-\rho_l)}$\, are chosen so
that
\begin{equation}\label{eq:rho:def}
  \xi_l(2\pi)=2\pi, \qquad\quad
  \frac{\min\{E_l(0),E_l(\pi)\}}{\max\{E_l(0),E_l(\pi)\}} = \rho_l.
\end{equation}
Mesh refinement near $x=0$ or $x=\pi$ causes the grid to sparsen near
$x=\pi$ or $x=0$ by a factor of
$\jd\max\{E_l(0),E_l(\pi)\} = 1+A_l\mu_l=
1/[1-\mu_l(1-\rho_l)]$ relative to a uniform grid.
Note that $\gridkap_l=0$ is a separate
case in (\ref{eq:xi:def}) and corresponds to uniform spacing.
For nonzero $\gridkap_l$, the mesh is uniform when $\rho_l=1$ and
becomes increasingly concentrated near the center ($\gridkap_l>0$) or
endpoints ($\gridkap_l<0$) as $\rho_l$ decreases toward 0.  Larger
values of $|\gridkap_l|$ lead to more localized regions of mesh
refinement.  In the present paper, we use either $\gridkap_l=0$ (uniform
  spacing), $\gridkap_l=2$ (mesh refinement near $x=\pi$), or
$\gridkap_l=-2$ (mesh refinement near $x=0,2\pi$).  Re-spacing the grid
from segment $l$ to $l+1$ amounts to interpolating the values of
$\eta$ and $\varphi$ from the old mesh to the new mesh.  We do this by
solving e.g.~$\eta\circ\xi_{l+1}(\alpha_j) = \eta\circ\xi_l(
  \xi_l^{-1}\circ\xi_{l+1}(\alpha_j))$ by Newton's method.

For time-stepping (\ref{eq:ww}), we use a 5th or 8th order Runge-Kutta
scheme \cite{hairer:I} in double-precision and a 15th order
spectral deferred correction method \cite{dutt, huang, minion,
  qadeer:faraday} in quadruple-precision (SDC15).  We also make use of
the 36th order filter popularized by Hou, Lowengrub and Shelley
\cite{hls94, hou97, HLS01, hou:li:07}.  This filter consists of
multiplying the $k$th Fourier mode by
\begin{equation}\label{eq:filter}
  \exp\left[-36\big(|k|/k_\text{max}\big)^{36}\right], \qquad
  k_\text{max} = M/2,
\end{equation}
which strikes a balance between suppressing aliasing errors
and resolving high-frequency modes. Here $M=M_l$ is the number of
grid points used on the $l^\text{th}$ spatial grid.

\subsection{Symmetry}
\label{sec:sym}

To compute symmetric standing waves \cite{penney:52, tadjbakhsh,
  concus:62}, we exploit a symmetry that reduces the
computation time by a factor of four \cite{mercer:92,water2}.  In
order for the solution to return to its initial state at time $T$, it
suffices to reach a rest state in which
\begin{equation}\label{eq:phi:T4}
  \varphi(x,T/4) = 0, \qquad 0\le x\le 2\pi,
\end{equation}
provided $\eta(x,0)$ and $\varphi(x,0)$ are spatially periodic even
functions satisfying
\begin{equation}\label{eq:shift}
  \eta(x+\pi,0) = \eta(x,0), \quad
  \varphi(x+\pi,0) = -\varphi(x,0).
\end{equation}
Indeed, reversing time about $t=T/4$ merely changes the sign of
$\varphi$, due to (\ref{eq:phi:T4}). Thus, $\eta(x,T/2)=\eta(x,0)$ and
$\varphi(x,T/2)=-\varphi(x,0)$, i.e.~the wave profile at $T/2$ is
identical to its initial configuration, while the velocity potential
changes sign.  By (\ref{eq:shift}), the evolution from $T/2$ to $T$
will be the same as that from $0$ to $T/2$ with a spatial translation
by $\pi$. Using (\ref{eq:shift}) again shows that
$\eta(x,T)=\eta(x,0)$ and $\varphi(x,T)=\varphi(x,0)$, i.e.~the
solution is time-periodic.

Once the initial conditions and period are found using symmetry to
accelerate the search for standing waves, we double-check that the
numerical solution evolved from $0$ to $T$ is indeed time-periodic.
Our Floquet analysis of the stability of these solutions is also
performed over the entire interval, $[0,T]$.

\subsection{Trust-region minimization}

Following Wilkening and Yu \cite{water2}, we compute time-periodic
solutions in a shooting framework by posing the problem as an
overdetermined nonlinear system of equations. The system is
overdetermined because we only solve for the leading Fourier modes
of the initial conditions, zero-padding the rest, but still impose
time-periodicity on all the Fourier modes. This improves robustness
over more traditional shooting methods and enables us to compute
sensitivities by discretizing linearized equations rather than
linearizing the discrete nonlinear equations. This allows us to
share the data structure for the Dirichlet-Neumann operator across
all columns of the Jacobian calculation, leading to a very efficient
parallel code framework.

Let $\hat\eta_k(t)$ and
$\hat\varphi_k(t)$ denote the Fourier modes of the solution.  For
symmetric standing waves, we assume the initial conditions have the
form
\begin{equation}\label{eq:init:stand}
  \begin{aligned}
    \hat\eta_k(0) &= c_{|k|}, \qquad (k=\pm2,\pm4,\pm6,\dots\;;\; |k|\le n), \\
    \hat\varphi_k(0) &= c_{|k|}, \qquad (k=\pm1,\pm3,\pm5,\dots\;;\; |k|\le n),
  \end{aligned}
\end{equation}
where $c_1,\dots,c_n$ are real numbers, and all other Fourier modes
are zero.  If the fluid depth is finite, we also set $\hat\eta_0=h$,
where $h$ is the mean fluid depth.  The integer $n$ controls the
cutoff beyond which the Fourier modes are zero-padded.  We normally
choose $n$ in the range $M/6\le n\le M/3$, where $M$ is the number of
grid points.  We also set
\begin{equation}\label{eq:c0:T}
  c_0=T,
\end{equation}
the unknown period.  We then define the objective function
\begin{equation}\label{eq:f:phi}
\begin{aligned}
  f(c) &= \frac{1}{2}r(c)^Tr(c) \approx \frac{1}{4\pi}\int_0^{2\pi}
  \varphi(x,T/4)^2\,dx, \\
  r_j &= \varphi(\alpha_j,T/4)/\sqrt{M}, \qquad
  \alpha_j = 2\pi j/M, \qquad
  0\le j<M.
\end{aligned}
\end{equation}
As explained in Section \ref{sec:sym} above, driving $f$ to zero will
yield time-periodic solutions with the symmetry properties expected of
standing waves.  Note that the number of equations, $m=M$, is
generally larger than the number of unknowns, $n$, due to zero-padding
of the initial conditions.  When computing extreme standing waves
using $f$ of the form (\ref{eq:f:phi}), we often refine the mesh as
$t$ approaches $T/4$, which further increases $m$ relative to $n$.
See Wilkening and Yu \cite{water2} for details.

We compute families of time-periodic solutions by specifying the mean
fluid depth $h$ (possibly infinite) and one of the parameters $c_k$
(with $k=k_0$, say) in a continuation algorithm. The other $c_k$ with
$k\ne k_0$ are chosen to minimize the objective function $f$ using a
variant of the Levenberg-Marquardt method \cite{water2}. The
Levenberg-Marquardt method requires a Jacobian matrix $J_{jk}=\pa
r_j/\pa c_k$, which we compute by solving the variational equations,
namely (\ref{eq:variational}) below, repeatedly, in parallel, with
different initial conditions; see Appendix~\ref{sec:J}.  If $f$
reaches a local minimum that is higher than a specified threshold
(normally $10^{-26}$ in double-precision and $10^{-52}$ in quadruple
  precision), we try again with a larger $n$, a larger $M$, a value of
$c_{k_0}$ closer to the last successful attempt, or a different index
$k_0$.  Switching the index is often useful when tracking a fold in
the bifurcation curve.  In some cases (see sections \ref{sec:shallow}
  and~\ref{sec:surf} below), rather than freeze the period or a
Fourier mode, we specify the value of $\eta(a,0)$ at $a=0$,
$a=\pi/2$ or $a=\pi$.  This is done by adding a component to $r$, namely
\begin{equation}\label{eq:eta:constr}
  r_m = \left[\hat\eta_0  +
    \left(\jt\sum_{k=1}^{M/2} \hat\eta_k(0)\big[2\cos(ka)\big]\right)
    - \eta(a,0)\right] w_m,
\end{equation}
where $\hat\eta_0$ and $\eta(a,0)$ are given, $w_m$ is a weight
(typically $10^{-2}$ or $1$), and $\hat\eta_k(0)$ is either 0 or one
of the $c_j$ (hence real).  In this approach, none of the $c_k$ are
frozen and $r(c)$ is a nonlinear function from $\mbb{R}^{n+1}$ to
$\mbb{R}^{m+1}$.

\section{Stability}
\label{sec:stable}

Once a time-periodic solution is found, we can check its stability
using Floquet theory \cite{coddington}. We restrict attention to
harmonic stability here, which is physically reasonable for symmetric
standing waves in a tank with vertical walls. Subharmonic stability
will be treated in a follow-up paper in the more general context of
traveling-standing water waves \cite{waterTS}, which include
traveling waves and standing waves as special cases.

Let $q(x,t)=\big(\eta(x,t);\varphi(x,t)\big)$ be a time-periodic
solution of the nonlinear system (\ref{eq:ww}) with period $T$. Here a
semicolon separates the components of a column vector. When
convenient, we drop the $x$-dependence in the notation and regard
$q(\cdot,t)$ (or simply $q(t)$) as evolving in an abstract Hilbert
space such as $L^2(\mbb T)\times H^1(\mbb T)$, where $\mbb T=\mbb
R\big/ 2\pi\mbb Z$. We write (\ref{eq:ww}) abstractly as
\begin{equation}\label{eq:qt}
  \partial_t q = F(q).
\end{equation}
The linearized system
\begin{equation}\label{eq:dot2}
  \partial_t \dot q = DF(q(t))\dot q
\end{equation}
generates a family of state-transition operators, $E$, such that
\begin{equation}
\dot q(t)=E(t,t_0)\dot q_0
\end{equation}
is the solution of (\ref{eq:dot2}) satisfying
$\dot q(t_0)=\dot q_0$.  Note that a dot represents a variational
derivative, not a time derivative; see (\ref{eq:eps:def}) in
Appendix~\ref{sec:J}. We denote the monodromy operator by
\begin{equation}\label{eq:ET:def}
E_T = E(T,0).
\end{equation}
%
It gives the change at time $T$ due to a perturbation at time 0.  For
$m\in\mathbb{Z}$ we have $E((m+1)T,mT)=E_T$ since $DF(q(t))$ is
periodic.  Formulas for the linearized system (\ref{eq:dot2})
are given by
\begin{equation}\label{eq:variational}
  \begin{aligned}
    \dot\eta_t &= \left(\dot\phi_y - \eta_x\dot\phi_x\right) - \left(\dot\eta \phi_x\right)', \\
    \dot\varphi_t &= P\left[ -\left(\dot\eta\phi_x\phi_y\right)' - \phi_x\dot\phi'
      + \phi_y\left(\dot\phi_y - \eta_x \dot\phi_x\right) - g\dot\eta +
      \frac\sigma\rho\pa_x\left( \frac{\dot\eta_x}{(1+\eta_x^2)^{3/2}}\right)\right],
  \end{aligned}
\end{equation}
where the right-hand side is evaluated at the free surface, subscripts
are partial derivatives taken before restricting to the boundary,
primes denote differentiation with respect to $x$ after restricting to
the boundary, e.g.
\begin{equation*}
  (\dot\eta\phi_x)' = \der{}{x}\big[ \dot\eta(x,t)\phi_x(x,\eta(x,t),t)\big] =
  \dot\eta_x\phi_x + \dot\eta\phi_{xx} + \dot\eta\phi_{xy}\eta_x,
\end{equation*}
and
$\dot\phi(x,y,t)=\der{}{\veps}\big\vert_{\veps=0}\phi(x,y,t;\veps)$,
with $\veps$ as in (\ref{eq:eps:def}) below. As a result,
$\dot\phi$ satisfies the Laplace equation with periodic boundary
conditions in $x$, $\dot\phi = \dot\varphi - \phi_y\dot\eta$ at the
free surface, and $\pa_y\dot\phi=0$ on the bottom boundary. Note that
$\dot\varphi$ is \emph{not} the restriction of $\dot\phi$ to the free
surface due to the $\phi_y\dot\eta$ term.  The linearized equations
(\ref{eq:dot2}) or (\ref{eq:variational}) are also solved when
computing the Jacobian in the Levenberg-Marquardt method; see
\cite{water2} for a detailed derivation of (\ref{eq:variational}) from
(\ref{eq:ww}) in this context.

The long-time behavior of solutions of (\ref{eq:dot2}) is governed by
the behavior of powers of $E_T$ due to the decomposition
\begin{equation}
  \dot q(t) = E(\tau,0)E_T^m\dot q_0, \qquad t = mT+\tau,
\end{equation}
where $m\in\mathbb{Z}$ and $0\le \tau<T$.
In finite dimensions, $E_T$ can be reduced to Jordan form, and its
eigenvalues are called Floquet multipliers.  In that case, all
solutions of (\ref{eq:dot2}) remain bounded iff (1) all the Floquet
multipliers satisfy $|\lambda_j|\le1$, and (2) only simple Jordan
blocks are associated with multipliers satisfying $|\lambda_j|=1$.
Moreover, $E(t,0)=P(t)e^{Rt}$, where $P(t)$ is periodic and $R=
\frac{1}{T}\log E_T$.
In infinite dimensions, a non-self-adjoint operator need not have any
eigenvalues at all.  Nevertheless, the growth of $E_T^l$ is connected
to the spectral radius via $\rho(E_T)=\lim_{m\rightarrow
  \infty}\|E_T^m\|^{1/m}$.  Moreover, if the point spectrum is
non-empty and lies partly outside the unit circle, or if a Jordan
chain of length greater than 1 exists for an eigenvalue on the unit
circle, we can construct unbounded solutions of (\ref{eq:dot2})
directly.

\subsection{Symmetry}

Since $E_T$ maps real-valued functions to real-valued functions, its
eigenvalues come in complex-conjugate pairs.  They also come in
inverse pairs when the time-periodic solution $q(t)$ has the
symmetry discussed in Section~\ref{sec:sym}.  To show this, we
denote the translation and Hamiltonian time-reversal operators by
\begin{equation}\label{eq:SR:def}
  S\left(\begin{array}{c} \eta(x,t) \\ \varphi(x,t) \end{array}\right)
   =\left(\begin{array}{c} \eta(x-\pi,t) \\ \varphi(x-\pi,t)\end{array}\right), \qquad
  R\left(\begin{array}{c} \eta(x,t) \\ \varphi(x,t) \end{array}\right)
   =\left(\begin{array}{c} \eta(x,-t) \\ -\varphi(x,-t) \end{array}\right).
\end{equation}
These operators commute with the evolution equations (\ref{eq:qt}) and
(\ref{eq:dot2}) in the sense that if $q$ and $\dot q$ are solutions, then
\begin{alignat*}{2}
  \partial_t (Sq)&=F(Sq), & \qquad \partial_t(S\dot q) &= DF(Sq(t))S\dot q, \\
  \partial_t (Rq)&=F(Rq), & \qquad \partial_t(R\dot q) &= DF(Rq(t))R\dot q.
\end{alignat*}
If $q=(\eta;\varphi)$ has initial conditions satisfying
(\ref{eq:shift}), then $Sq=Rq$ and $q=SRq$.  It follows that $SR\dot
q$ satisfies (\ref{eq:dot2}) provided that $\dot q$ does.  Thus, if
$\dot q(0)$ is an eigenfunction of $E_T$, i.e.
\begin{equation}
  \dot q(T) = \lambda \dot q(0),
\end{equation}
then $\dot q(-T)=\lambda^{-1}\dot q(0)$ and
$SR\dot q(0)$ is an eigenfunction with eigenvalue $\lambda^{-1}$:
\begin{equation}
  E_T[SR\dot q(0)] =
  \left(\begin{array}{c} \dot\eta(x-\pi,-T) \\ -\dot\varphi(x-\pi,-T)
    \end{array}\right) =
  \lambda^{-1}\left(\begin{array}{c} \dot\eta(x-\pi,0) \\ -\dot\varphi(x-\pi,0)
  \end{array}\right) =
  \lambda^{-1}SR\dot q(0).
\end{equation}
When a Floquet multiplier $\lambda$ is on the unit circle,
$\lambda^{-1}=\overline\lambda$. Thus, if the eigenspace of $\lambda$
is one-dimensional, that of $\overline\lambda$ will also be
one-dimensional and $\overline{\dot q(0)}$ will be a complex multiple
of $SR\dot q(0)$.  Floquet multipliers not on the unit circle come in
inverse pairs of real numbers, or in groups of four complex numbers:
$\lambda$, $\lambda^{-1}$, $\overline\lambda$,
$\overline{\lambda^{-1}}$.

So far we have only made use of property
(\ref{eq:shift}). An additional symmetry can
be exploited to untangle even eigenfunctions from odd ones.
Inspecting (\ref{eq:variational}), we see that if $\eta_0$ and
$\varphi_0$ are both even functions while $\dot\eta_0$ and
$\dot\varphi_0$ are both even or both odd, then $\dot\eta$ and
$\dot\varphi$ will remain even or odd for all time.  As a result,
$E_T$ leaves invariant the subspace of pairs
$(\dot\eta_0;\dot\varphi_0)$ that are both even or both odd functions.
Thus, we can
search for even and odd eigenfunctions of $E_T$ separately, rather
than all at once.  This leads to a more canonical set of
eigenfunctions as well as a fourfold computational savings, as it is
less expensive to diagonalize two matrices of half the size.

\subsection{A Fourier basis}
\label{sec:fourier}

One way to approximate $E_T$ is to evolve (\ref{eq:variational}) with
$\dot\eta_0$ and $\dot\varphi_0$ ranging over all discrete
$\delta$-functions defined on the mesh.  This is the approach taken by
Mercer and Roberts \cite{mercer:92}. The discrete $\delta$-functions
can be thought of as ``trigonometric Lagrange polynomials,'' so as
long as $E_T$ acts on interpolated values of a smooth function, the
method can be quite accurate. Nevertheless, we find that it is more
efficient and robust to express $E_T$ in a Fourier basis before
computing eigenvalues.  To keep the matrix entries real and avoid
storing negative index Fourier modes, we use the r2c basis of FFTW,
with basis functions
\begin{equation}
  1\;,\; 2\cos x\;,\; -2\sin x\;,\; \dots\;,\;
  2\cos kx\;,\; -2\sin kx\;,\; \dots.
\end{equation}
Since the $k=0$ modes of $\dot\eta(x,t)$ and $\dot\varphi(x,t)$ remain
constant in time in (\ref{eq:variational}) and the right-hand side of
(\ref{eq:variational}) is unchanged by adding constants to $\dot\eta$
and $\dot\varphi$,
$E_T = E_T^0\oplus E_T^\#$, where $E_T^0$ is the identity map on the
two-dimensional space of constant functions
$(\dot\eta_\text{const};\dot\varphi_\text{const})$, and $E_T^\#$ is
the restriction of $E_T$ to functions $(\dot\eta(x)$;$\dot\varphi(x))$
in which both components have zero mean.  Next we define $\hat E_T =
QE_T^\#Q^{-1}$, where
\begin{equation} \label{eq:Q:def}
  \begin{aligned}
    &[Q^{-1}z](x) = 2\sum_{k=1}^\infty \left(\begin{array}{c}
       z_{4k-3}\cos kx - z_{4k-2}\sin kx \\[2pt]
       z_{4k-1}\cos kx - z_{4k}\sin kx\end{array}
       \right), \quad
    \begin{array}{rl}
      \longleftarrow & \dot\eta(x) \\ \longleftarrow & \dot \varphi(x)
    \end{array}\\
    &[Q\dot q]_{4j+j'} = \left\{\begin{array}{cl}
      \frac{1}{2}\big[\wdg{\dot\eta}_j+\wdg{\dot\eta}_{-j}\big],
      & j'=-3 \\[3pt]
      \frac{1}{2i}\big[\wdg{\dot\eta}_j-\wdg{\dot\eta}_{-j}\big],
      & j'=-2 \\[3pt]
      \frac{1}{2}\big[\wdg{\dot\varphi}_j+\wdg{\dot\varphi}_{-j}\big],
      & j'=-1 \\[3pt]
      \frac{1}{2i}\big[\wdg{\dot\varphi}_j-\wdg{\dot\varphi}_{-j}\big],
      & j'=0 \\[3pt]
    \end{array}\right\}, \qquad (j\ge1).
  \end{aligned}
\end{equation}
When applied to a real-valued pair of functions
$(\dot\eta;\dot\varphi)$, considered as functions of $x$ only, $Q$ simply
extracts the real and imaginary parts of their Fourier modes and
interlaces them. To compute column $4k+k'$ of $\hat E_T$, we evolve
(\ref{eq:variational}) with initial conditions
\begin{equation}\label{eq:q0:dot}
  \dot q_0 = \left\{\begin{array}{rl}
      \big(2\cos kx;0\big), & k'=-3 \\[2pt]
      \big(\text{$-2$}\sin kx;0\big), & k'=-2 \\[2pt]
      \big(0;2\cos kx\big), & k'=-1 \\[2pt]
      \big(0;-2\sin kx\big), & k'=0
    \end{array}\right\}
\end{equation}
to time $T$ and extract real and imaginary parts of
$\wdg{\dot\eta}_j(T)$, $\wdg{\dot\varphi}_j(T)$ to obtain rows
$4j-3,\dots,4j$.  Note that $Q$, $E_T^\#$ and $\hat E_T$ act on the
real and imaginary parts of complex data separately. It was important
to avoid complex conjugates in (\ref{eq:Q:def}) to ensure that $Q$ is
complex linear.


\subsection{Stability of the zero-amplitude solution}\label{sec:zro:stab}

When linearizing about the flat rest state, we can compute $\hat E_T$
explicitly. The linearized equations are
\begin{equation}\label{eq:lin}
 \partial_t\dot\eta = \mc{G}\dot\varphi, \qquad
 \partial_t\dot\varphi = P[-g\dot\eta],
\end{equation}
where the Dirichlet-Neumann operator is given by
\begin{equation}
  \mc{G}\big[e^{ikx}\big] = \big[k\tanh k\hh\big]e^{ikx}.
\end{equation}
The four initial conditions listed in (\ref{eq:q0:dot}) lead to
the solutions
\begin{equation*}
  \begin{aligned}
    \dot\eta(x,t) = \cos(\omega_k t-\theta_{k'})\left[ \gamma_{kk'} e^{ikx} +
      \bar\gamma_{kk'} e^{-ikx}\right], \\
    \dot\varphi(x,t) = -(g/\omega_k)\sin(\omega_k t-\theta_{k'})
    \left[ \gamma_{kk'} e^{ikx} +
      \bar{\gamma}_{kk'} e^{-ikx}\right],
  \end{aligned}
\end{equation*}
where $\omega_k=\sqrt{kg\tanh kh}$ and $\gamma_{kk'}$, $\theta_{k'}$ are given by
\begin{equation}
  \begin{array}{c||c|c|c|c} k' & \;-3\; & \;-2\; & -1 & 0 \\ \hline
    \gamma_{kk'} & 1 & i & \;\omega_k/g\; & \;i\omega_k/g \\
    \theta_{k'} & 0 & 0 & \pi/2 & \pi/2.
  \end{array}
\end{equation}
It follows that $\hat E_T$ is block-diagonal with $4\times 4$ blocks
given by
\begin{equation*}
  F_k =
  \left(\begin{array}{c@{}c@{}c@{}c}
    \cos\omega_k T & 0 & \frac{\omega_k}{g}\sin\omega_k T & 0 \\
    0 & \cos\omega_k T & 0 & \frac{\omega_k}{g}\sin\omega_k T \\
    -\frac{g}{\omega_k} \sin\omega_k T & 0 & \cos\omega_k T & 0 \\
    0 & -\frac{g}{\omega_k} \sin\omega_k T & 0 & \cos\omega_k T
    \end{array}\right), \qquad (k\ge1).
\end{equation*}
A natural choice of period is $T=2\pi/\omega_1$, which amounts to
linearizing about an infinitesimally small standing wave of wavelength
$2\pi$.  The $k$th block of $\hat E_T$ has two double-eigenvalues (or
  one quadruple-eigenvalue) $\lambda_j = e^{\pm i\omega_kT}$, where
$j\in\{4k-3,\dots,4k\}$ enumerates the eigenvalues (counting
  multiplicity).  Indeed, $F_k = U_k \Lambda_k U_k^{-1}$ with
\begin{equation*}
  U_k = \begin{pmatrix}
    i & -i & 0 & 0 \\
    0 & 0 & i & -i \\
    \beta_k & \beta_k & 0 & 0 \\
    0 & 0 & \beta_k & \beta_k
  \end{pmatrix}, \qquad \Lambda_k = \opn{diag}
  \begin{pmatrix}
    c_k-is_k \\ c_k+is_k \\
    c_k-is_k \\ c_k+is_k \end{pmatrix},
\end{equation*}
where $c_k=\cos \omega_kT$, $s_k=\sin\omega_kT$, and
$\beta_k=g/\omega_k$.  The infinite matrix $U$ containing $4\times 4$
diagonal blocks $U_k$ is invertible as an operator from
$l^2(\mathbb{N})$ to $l^2(\mathbb{N},\mu)$, with measure $\mu$ chosen
to obtain the weighted norm
\begin{equation}\label{eq:z:mu}
  \|z\|_\mu^2 = \sum_{k=1}^\infty \bigg[
    \Big(|z_{4k-3}|^2 + |z_{4k-2}|^2 \Big)
    + (\omega_k/g)^2 \Big(|z_{4k-1}|^2 + |z_{4k}|^2\Big)\bigg].
\end{equation}
We note that $Q$ is an isomorphism from $L^2_\#(\mathbb{T})\times
H^1_\#(\mathbb{T})$ to $l^2(\mathbb{N},\mu)$, where
$L^2_\#(\mathbb{T})$ is the space of square-integrable functions with
zero mean, and $H^1_\#(\mathbb{T})$ is the Sobolev space of
$2\pi$-periodic functions with zero mean and a square-integrable weak
derivative.  Thus, when linearizing about the zero-amplitude solution,
$E_T^\#$ is an isomorphism on $L^2_\#(\mathbb{T})\times
H^1_\#(\mathbb{T})$.  Since all the eigenvalues lie on the unit
circle, the zero-amplitude solution is linearly stable to harmonic
perturbations both forward and backward in time.

\subsection{Numerical procedure for computing Floquet multipliers}

To determine the stability of finite-amplitude solutions, the Floquet
multipliers must be computed numerically. This is more difficult for
standing waves than for traveling waves because the eigenvalues of
the monodromy operator $E_T$ are complex numbers of the form
\begin{equation}\label{eq:sig:def}
  \lambda_j = |\lambda_j|e^{i\sigma_j},
\end{equation}
but the eigenfrequencies $\sigma_j$ are only known modulo $2\pi$.  By
contrast, for traveling waves, these eigenfrequencies can be computed
directly, giving a convenient ordering of the eigenvalues and a means
of extracting the ``leading'' eigenvalues in finite dimensional
approximations of the eigenvalue problem.  Our solution is to define a
mean wave number of the corresponding eigenfunction that replaces the
eigenfrequency in ordering the eigenvalues. Optionally, as a second
step, we re-order the eigenvalues by performing a numerical
continuation back to the zero-amplitude solution of
Section~\ref{sec:zro:stab} above. This involves developing an
algorithm to match eigenvalues of monodromy operators of nearby
time-periodic solutions.  Before attempting this matching, we need to
ensure that most of the eigenvalues we have computed are well
resolved. We find that working in Fourier space is less expensive than
the discrete $\delta$-function approach \cite{mercer:92}, and provides
a natural means of discarding spurious eigenvalues.

\begin{figure}[t]
  \begin{algorithm} \label{alg:fmult}
    (Floquet multipliers)\\
    \begin{itemize}[leftmargin=24pt,rightmargin=24pt]
  \item[] Goal: compute the leading $n^*$ eigenvalues of $E_T$, ordered
    by mean wave number (\ref{eq:mwn}) or by homotopy to the flat rest state.
  \item[(0)] Choose an initial $M$ and $n$. Here $M$ is the number of
    grid points, $n=4k_\text{max}$ is the number of different
    initial conditions $\dot q_0$ in (\ref{eq:q0:dot}) that will be
    evolved, and $k_\text{max}$ is a Fourier cutoff.  (Typically
      $n\approx 2n^*$, $M\approx 3n/2$.)
  \item[(1)] Compute the leading $(2M-4)\times n$ sub-block $J$ of
    $\hat{E}_T$ by numerically solving (\ref{eq:variational}) with initial
    conditions (\ref{eq:q0:dot}).  Here $2M-4$ is the number of
    degrees of freedom in $\dot q=(\dot\eta,\dot\varphi)$, excluding the
    constant and Nyquist modes. As explained in Section~\ref{sec:fourier},
    the matrix entries of $\hat E_T$, and hence $J$, are real.
  \item[(2)] Compute the eigenvalues and eigenvectors of
    $J_{1:n,1:n}$. We modified the Template Numerical Toolkit eigenvalue
    solver to work in quadruple-precision, though most of the
    computations of this paper are done in double-precision using
    LAPACK.
  \item[(3)] Sort the eigenvalues $\lambda_j$ and eigenfunctions
    $z^\e j$ by mean wave number (\ref{eq:mwn}).  Also compute the
    residual error
  \begin{equation}\label{eq:err}
    err_j = \big\|Jz^\e j -\lambda_j[z^\e j;0]\big\|,
  \end{equation}
  where $[z^\e j;0]$ extends $z^\e j$ from $\mathbb{C}^n$ to $\mathbb{C}^{2M-4}$ by
  zero-padding.
\item[(4)] If necessary, repeat steps (1)--(4) on a finer mesh (increasing $M$ and
  $n$) until the first $n^*$ errors are smaller than a user specified
  error tolerance.
\item[(5)] (optional) If a family of monodromy operators and their Floquet
  multipliers have been computed, use the matching algorithm described
  in Appendix~\ref{sec:matching} to track individual eigenvalues as
  the parameter changes. This involves some re-ordering of eigenvalues,
  and can cause a few exchanges across the boundary between leading
  and discarded eigenvalues.  
    \end{itemize}
    \begin{equation*}
      \rule{5in}{.5pt}
    \end{equation*}
  \end{algorithm}
\end{figure}

The new algorithm is summarized in Algorithm~\ref{alg:fmult}.  In step
(1), the computational savings over the $\delta$-function approach
comes from choosing $n/(2M)\approx 1/3$ rather than $n=2M$. In other
words, we evolve only those initial conditions that are smooth enough
(relative to the mesh spacing) to yield accurate answers. It is also
useful to compute $J$ from step (1) in batches of columns by solving
(\ref{eq:variational}) in parallel with multiple right-hand sides, as
explained in Appendix~\ref{sec:J}.  In step (3), we sort the
eigenvalues of $\hat E_T$ by mean wave number, $\langle k\rangle$,
which measures how oscillatory the corresponding eigenfunctions are.
More specifically, if $z=z^\e{j}$ is an eigenvector of $\hat E_T$
corresponding to the eigenfunction $\dot
q=(\dot\eta;\dot\varphi)=Q^{-1}z^\e j$ of $E_T$, we define $\langle
k\rangle =\langle k\rangle_j$ by
\begin{align}
\notag
  \langle k\rangle &= \frac{
  \sum_{k=1}^\infty k\left[
    \left|\wdg{\dot \eta}_k\right|^2 +
    \left|\wdg{\dot \eta}_{-k}\right|^2 +
    \frac{\omega_k^2}{g^2}\left(
    \left|\wdg{\dot \varphi}_k\right|^2 +
    \left|\wdg{\dot \varphi}_{-k}\right|^2
  \right)\right] }{
  \sum_{k=1}^\infty \left[
    \left|\wdg{\dot \eta}_k\right|^2 +
    \left|\wdg{\dot \eta}_{-k}\right|^2 +
    \frac{\omega_k^2}{g^2}\left(
    \left|\wdg{\dot \varphi}_k\right|^2 +
    \left|\wdg{\dot \varphi}_{-k}\right|^2
    \right)\right]} \\
  \label{eq:mwn}
  &= \frac{1}{\|z\|_\mu^2} \sum_{k=1}^\infty k\bigg[
  \Big(|z_{4k-3}|^2 + |z_{4k-2}|^2\Big) +
      (\omega_k/g)^2 \Big(|z_{4k-1}|^2 + |z_{4k}|^2\Big)\bigg],
\end{align}
where $\|z\|_\mu$ was defined in (\ref{eq:z:mu}).  This gives an
unambiguous order to the eigenvalues, where $\lambda_j$ precedes
$\lambda_l$ if $\langle k\rangle_j<\langle k\rangle_l$, that roughly
coincides with how difficult they are to compute --- more oscillatory
eigenfunctions require additional resolution.  Ties do occur, but
generally in clusters of two or four that should be grouped together
anyway to form a single invariant subspace as they correspond to
complex conjugate pairs and time-reversed dynamics.  For example, in
the zero-amplitude case discussed above, $\langle k\rangle=k_0$ for
all four eigenvectors associated with block $F_{k_0}$. The same result
would be obtained by ordering the eigenvalues $\lambda_j=e^{\pm
  i\omega_{k_0}T}$, $j\in\{4k_0-3,\dots,4k_0\}$, by eigenfrequency
$\omega_{k_0}$; however, these are only known modulo $2\pi$ in the
numerical method without also looking at the eigenvectors.

We will see below that the residual errors grow as the eigenfunctions
become more oscillatory and shrink as $n$ and $M$ are increased.
Thus, by ordering the eigenvectors by mean wave number, the more
accurately computed eigenvalues will appear at the front of the list.
The residual error of a fixed number of them, $n^*$, can be made
smaller than a specified tolerance by increasing $n$ and $M$
sufficiently.

\section{Results}
\label{sec:results}

\subsection{Stability of standing waves on deep water}
\label{sec:deep}

\begin{figure}
\begin{center}
\includegraphics[width=.6\linewidth]{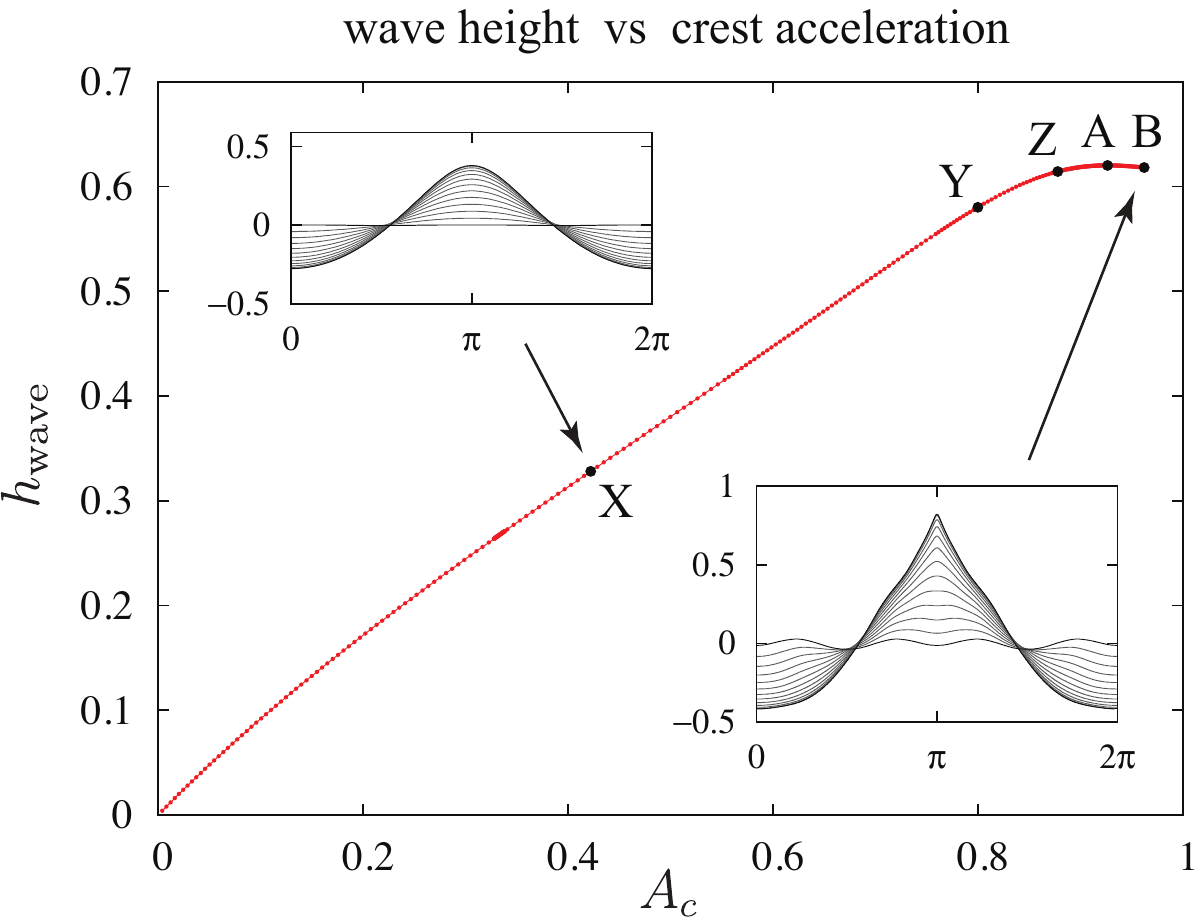}
\end{center}
\caption{\label{fig:bif:00} Wave height $h_\text{wave}$ versus crest acceleration
  $A_c$ over the range $0\le A_c\le 0.9639$, along with snapshots of
  $\eta(x,t)$ for solutions X and B over a quarter-period. For small
  values of $A_c$, the solution is nearly sinusoidal ($\eta\approx
  -A_c\cos x\sin t$, $\varphi\approx -A_c\cos x \cos t$, $h_\text{wave}\approx
  A_c$).  As $A_c$ increases, the wave crest sharpens and the surface
  dynamics become more complicated.  }
\end{figure}

In this section we present a Floquet analysis of standing waves in
deep water over the range $0\le A_c \le 0.9639$, where $A_c$ is the
(downward) acceleration of a fluid particle at the wave crest when it
reaches maximum height, normalized so that the acceleration of gravity
is $g=1$.  A bifurcation curve showing wave height $h_\text{wave}$
versus crest acceleration
for these standing waves is given in Fig.~\ref{fig:bif:00}.  The wave
height, defined here as half the vertical crest-to-trough distance,
reaches a local maximum of $0.62017$ at solution A, where
$A_c=0.92631$.  The other solutions were selected arbitrarily to
represent typical solutions along the bifurcation curve, and are
labeled cyclically to match high-resolution plots of solutions A and B
in previous work \cite{water1}.  Beyond solution B, for values of
crest acceleration in the range $0.98<A_c<1$, the bifurcation curve
breaks up into several disjoint branches and ceases to be in 1-1
correspondence with crest acceleration \cite{water1}.  Rather than
sharpen to a corner or cusp as conjectured by Penney and Price
\cite{penney:52, grant, okamura:98}, the solution develops fine-scale
oscillatory structure near the crest tip and throughout the domain
\cite{water1,water2}.

\begin{figure}
\begin{center}
\includegraphics[width=.8\linewidth]{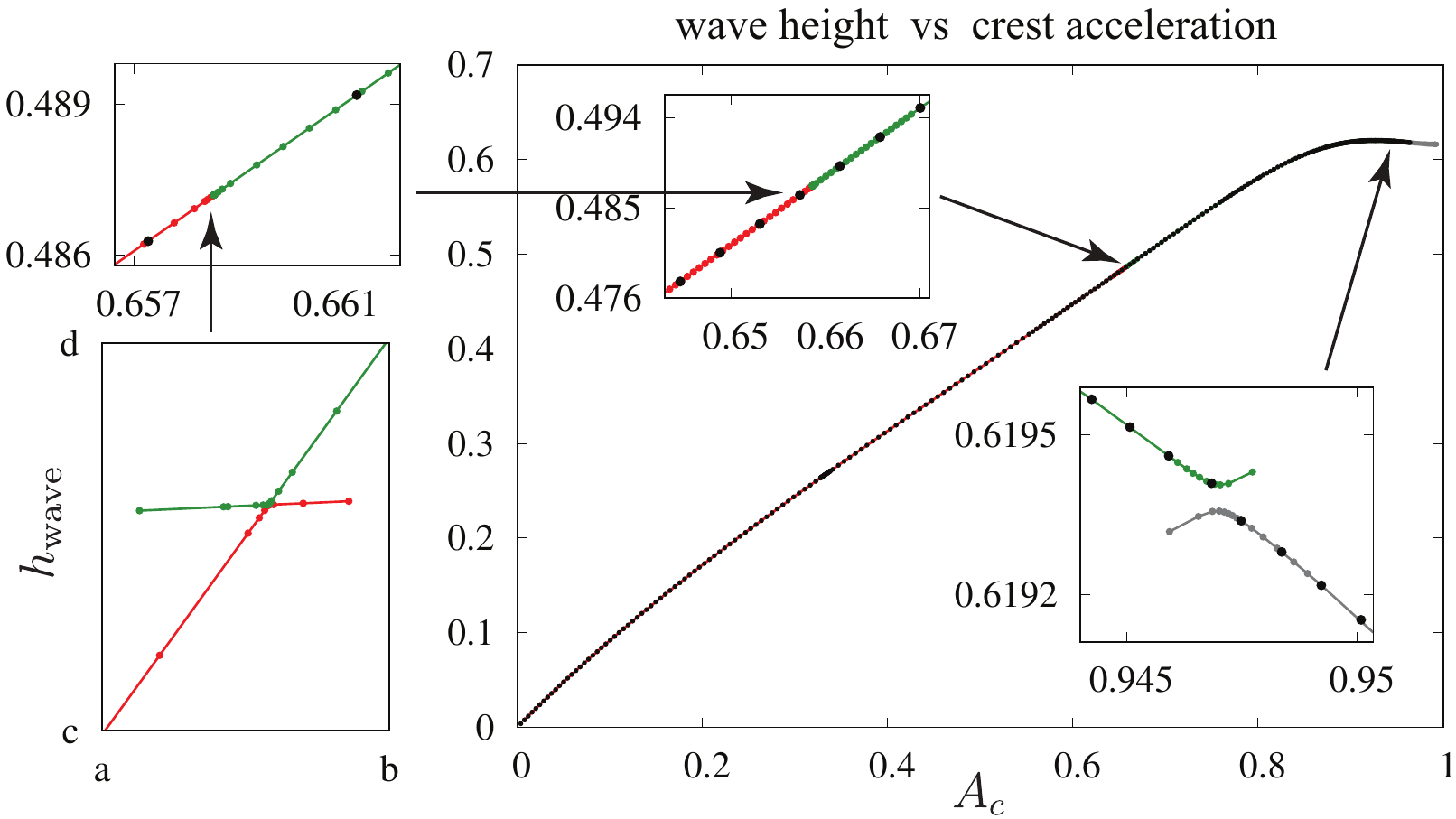}
\end{center}
\caption{\label{fig:bif:00:jump} Two ``resonant disconnections'' in
  the bifurcation diagram over the range $0<A_c<0.9639$ were found in
  \cite{water2}. The first one, with $A_c$ between
  $a=0.658\,621\,068\,280$ and $b=0.658\,621\,070\,350$ and
  $h_\text{wave}$ between $c=0.487\,198\,766\,380$ and
  $d=0.487\,198\,767\,760$, was only visible when the bifurcation
  curve was computed in quadruple-precision. The black markers denote
  points where Floquet multipliers are computed in the present study,
  while the red, green and gray markers denote the solutions computed
  in \cite{water2} with no stability calculations performed. }
\end{figure}

Our interest here is in the stability of the solutions at lower
amplitude, where the bifurcation curve is either smooth or contains
resonant disconnections that are too weak to be observed in numerical
simulations to date. In a careful search for such disconnections in
quadruple precision, Wilkening and Yu \cite{water2} were only able to
find two in the range $0\le A_c\le 0.9639$, namely at $A_c=0.658621$
and $A_c=0.947$. We note that to see the resonant disconnection at
$A_c=0.658621$ in Figure~\ref{fig:bif:00:jump}, we had to zoom in to a
window $[a,b]\times[c,d]$ with $\frac ba-1\approx3.1\times10^{-9}$ and
$\frac dc-1\approx 2.8\times10^{-9}$. Thus, as noted in \cite{water2},
it is extremely unlikely that this resonance could be computed in
double-precision without knowing where to look.  Since the effect of
such resonances is highly localized and our continuation path jumps
over the two disconnections we know about without refining the step
length to resolve them, we do not observe any unusual behavior in the
dependence of the multipliers on crest acceleration at $A_c=0.658621$
or at $A_c=0.947$ in Figs.~\ref{fig:mag:00}--\ref{fig:arg:00} below.
It is beyond the scope of the present work to quantify the effect of
these disconnections on the Floquet multipliers. We emphasize that the
stability question is separate from the question of whether solutions
occur in smooth families \cite{plotnikov01, iooss05}. Given a
time-periodic solution, steps (1)--(4) of Algorithm~\ref{alg:fmult}
will compute its leading Floquet multipliers.  Only in the optional
step (5) do we assume that nearby parameter values should lead to
small changes in the eigenvalues.

\begin{table}[t]
  \caption{\label{tab:params}
    Parameters used to compute Floquet
    multipliers of standing water waves.}
\begin{equation*}
  \begin{array}{c|c|c|c|c|c|c}
    \text{$A_c$ range} & M & N & n & n^* & \text{scheme} & \text{\#} \\ \hline
    \text{0--0.8871}      & 768  & 1440 &  600 & 300 &  5d & 297 \\
    \text{0.8878--0.9295} & 1024 & 1920 &  600 & 300 &  5d & 41 \\
    \text{0.9301--0.9639} & 1536 & 2880 &  900 & 300 &  5d & 42 \\
    0.8776355\,(Z)        & 1024 &  384 &  800 & 400 & 15q & 1 \\
    0.9263124\,(A)        & 1536 & 1200 & 1000 & 500 & 15q & 1 \\
    0.9620239\,(B)        & 3840 & 1152 & 1500 & 800 &  8d & 1
  \end{array}
\end{equation*}
\end{table}

The parameters used to generate the matrices $\hat E_T$ for this
stability study are given in Table~\ref{tab:params}.  The schemes $5d$
and $8d$ in the table correspond to Dormand and Prince's 5th and 8th
order Runge-Kutta methods \cite{hairer:I}, respectively.  (The $d$
stands for double-precision).  $15q$ is a 15th order spectral-deferred
correction method \cite{minion,qadeer:faraday}, with computations
performed in quadruple-precision.  The column labeled `\#' gives the
number of solutions that were computed in the given range of crest
acceleration values. The points are uniformly distributed in each
range, except near $A_c=0.333$, where additional points were used to
resolve a bubble of instability that we observed there.  The last 3
entries in the table are re-computations of solutions Z, A, and B on a
finer grid to check that the Floquet multipliers do not change when
the mesh is refined. For solutions Z and A, these re-computations were
done in quadruple precision.

Figures \ref{fig:mag:00}--\ref{fig:arg:00} show the first 300 Floquet
multipliers $\lambda_j = |\lambda_j|e^{i\sigma_j}$ for each of the 380
solutions in the first 3 rows of Table~\ref{tab:params}.  We actually
used $n^*=360$ for steps (1)--(4) of Algorithm~\ref{alg:fmult}, which
gives 360 eigenvalues sorted by mean wave number. We then re-ordered
these using the matching algorithm described in
Appendix~\ref{sec:matching} to obtain continuous curves depending on
the parameter $A_c$. Once matched, we kept the first 300, which is
sometimes a slightly different set than the 300 eigenvalues of smallest
mean wave number.  Multipliers not on the unit circle correspond to
unstable modes either forward or backward in time. They are plotted
with black markers in these figures.  Solutions X, Y, Z, A and B are
plotted with green markers.


\begin{figure}
\begin{center}
\includegraphics[width=.6\linewidth]{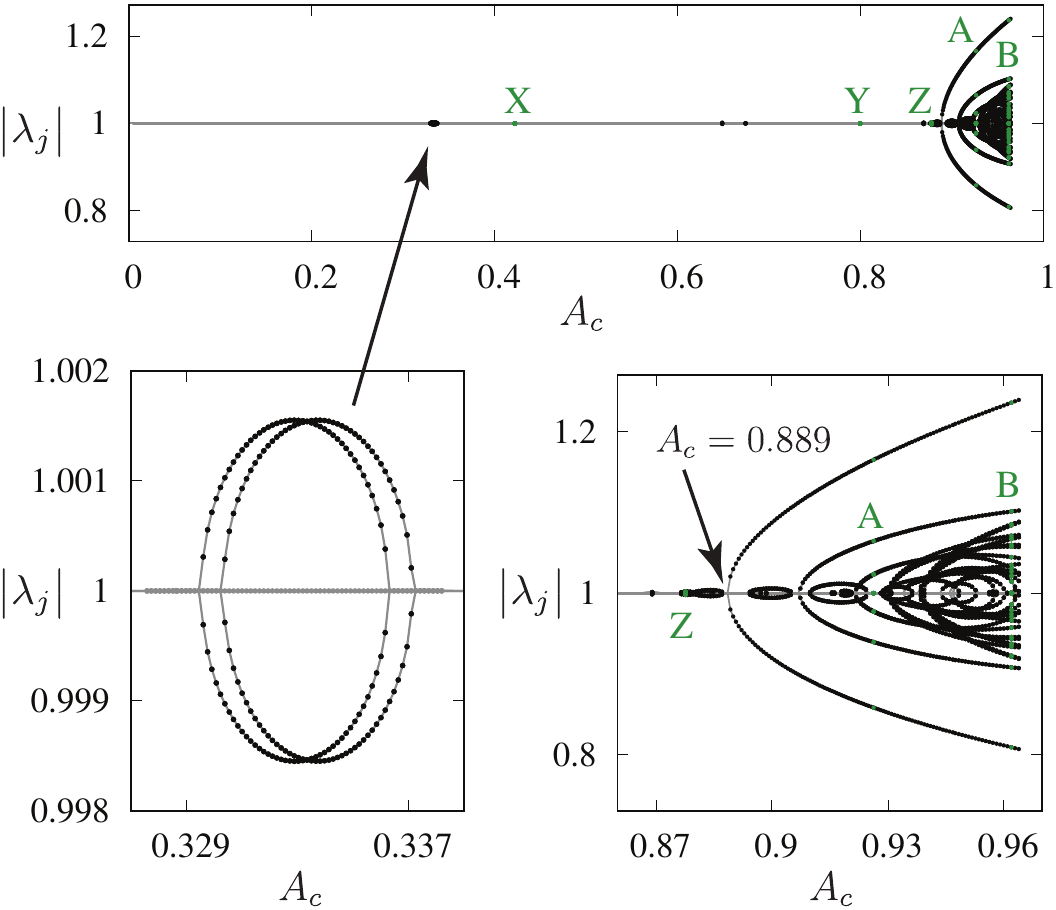}
\end{center}
\caption{\label{fig:mag:00} Magnitudes of the first 300 Floquet
  multipliers of the 380 solutions listed in rows 1--3 of
  Table~\ref{tab:params}.  We found several bubbles of instability
  below the first major unstable branch observed by Mercer and Roberts
  \cite{mercer:92} at $A_c\approx 0.889$.}
\end{figure}

\begin{figure}
\begin{center}
\includegraphics[width=.8\linewidth]{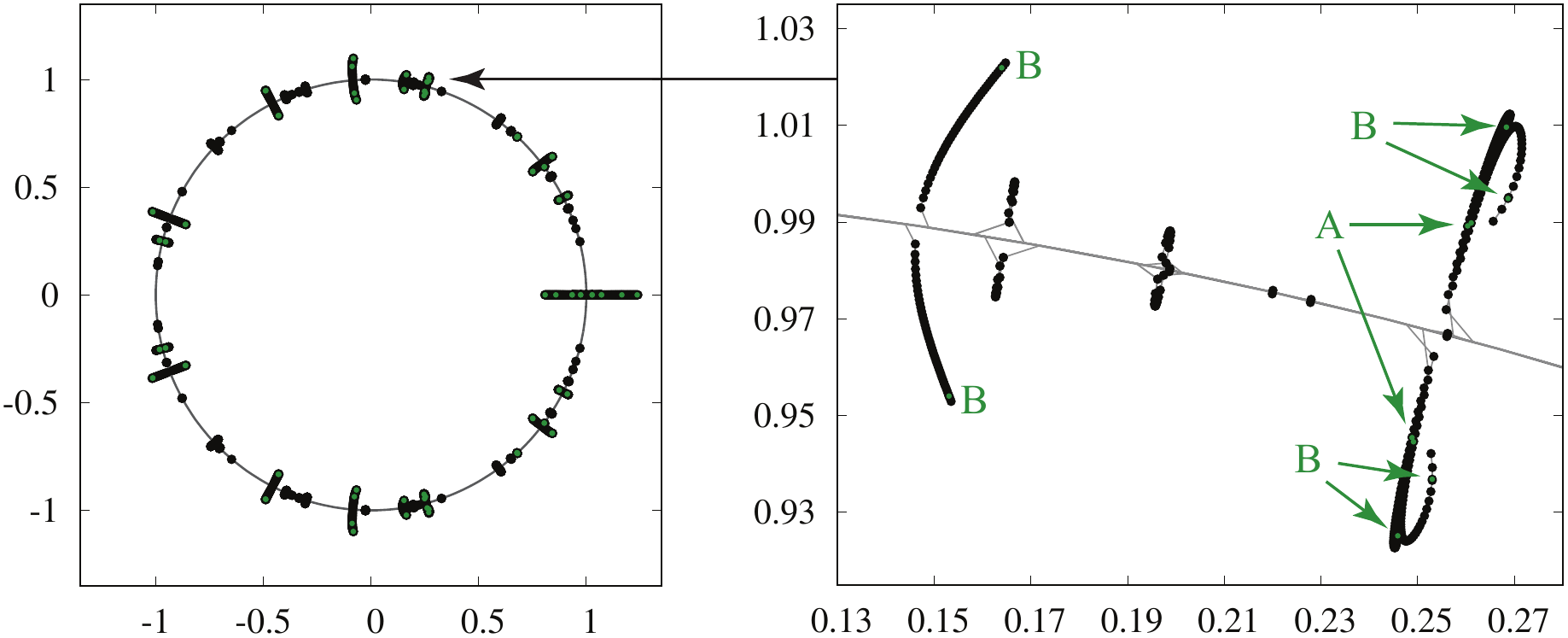}
\end{center}
\caption{\label{fig:mag:00a} Position in the complex plane of the
  Floquet multipliers of Fig.~\ref{fig:mag:00}. Black markers denote
  multipliers whose distance to the unit circle exceeds $10^{-5}$.
  All values of $A_c$ are plotted simultaneously. }
\end{figure}


\begin{figure}
\begin{center}
\includegraphics[width=.67\linewidth]{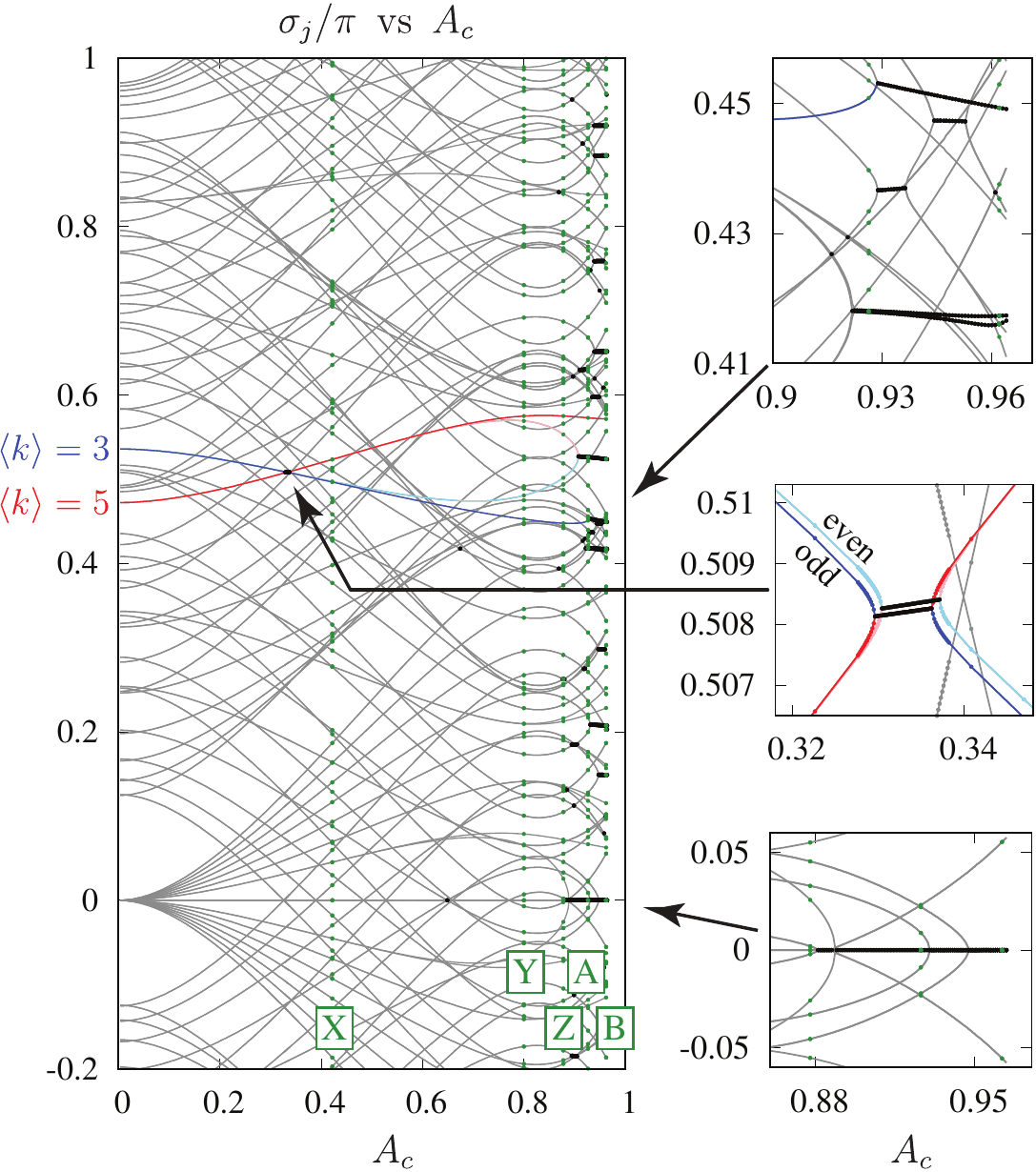}
\end{center}
\caption{\label{fig:arg:00} Dependence of the first 300
  eigenfrequencies on crest acceleration.  Since Floquet multipliers
  come in conjugate pairs, it was not necessary to plot the entire
  interval $-1\le\sigma/\pi\le 1$.  Unstable modes are generated when
  two eigenfrequencies coalesce and the corresponding multipliers
  leave the unit circle (see Figs.~\ref{fig:mag:00}
  and~\ref{fig:mag:00a}). These pairs often re-gain stability by
  joining the unit circle elsewhere, which allows the eigenfrequencies
  to split apart again.  }
\end{figure}

\begin{figure}
\begin{center}
\includegraphics[width=.7\linewidth]{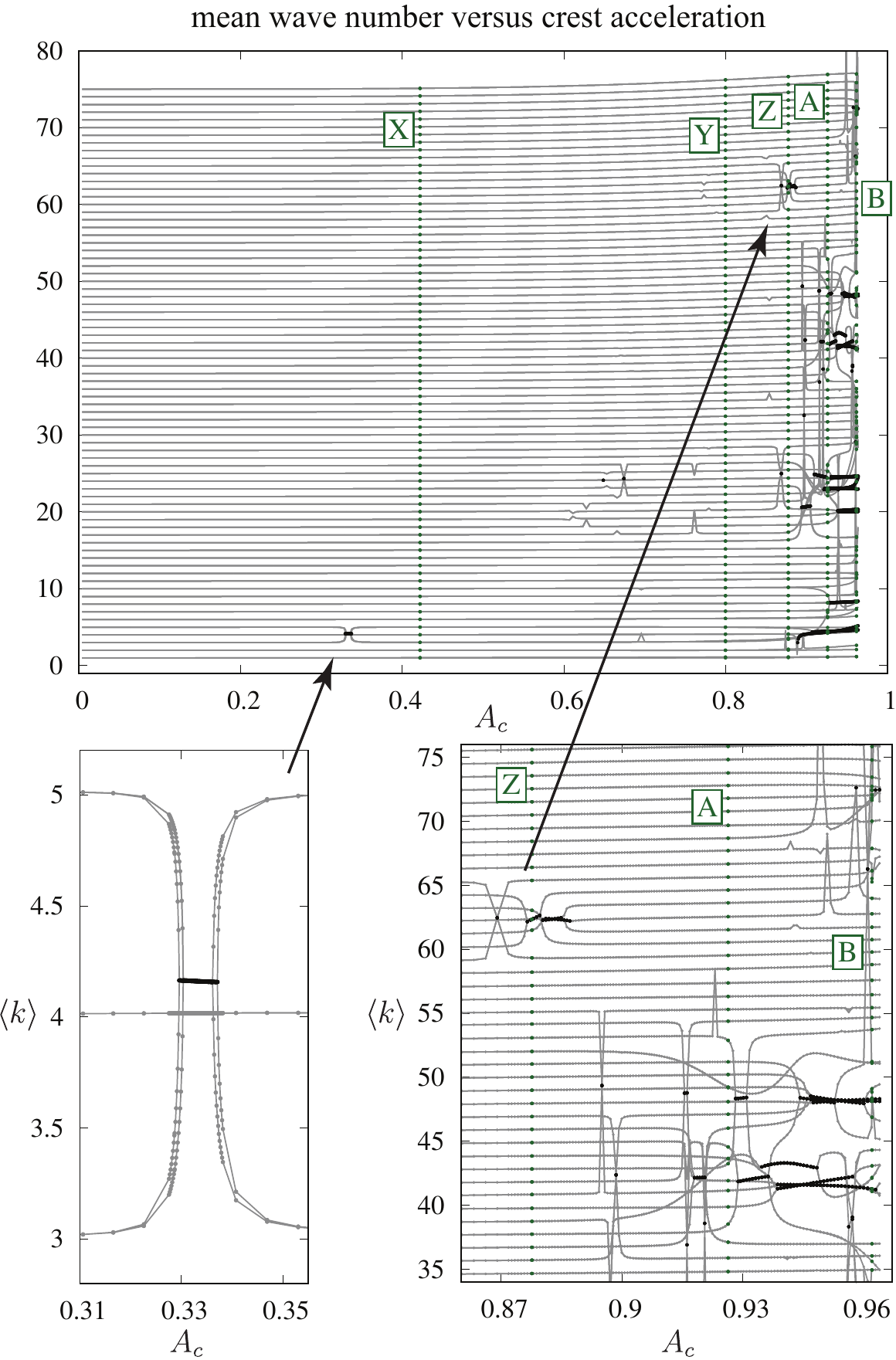}
\end{center}
\caption{\label{fig:kk:00}
  Mean wave number versus crest acceleration for the first 300
  Floquet multipliers. The matching algorithm of Appendix~\ref{sec:matching}
  was used to connect the eigenvalues of adjacent problems into
  continuous curves. Bubbles of instability form when eigenvalues
  collide and leave the unit circle. Taken together, plots of mean
  wave number, eigenfrequency (Fig.~\ref{fig:arg:00}), and magnitude
  (Fig.~\ref{fig:mag:00}) give a complete picture of these nucleation
  events.}
\end{figure}

Figure~\ref{fig:mag:00} shows the magnitude of the Floquet multipliers
as a function of $A_c$. Mercer and Roberts \cite{mercer:92} found that
standing waves with crest acceleration below $A_c=0.889$ are stable
to harmonic perturbation. Our results agree with theirs in that we
observe a dominant branch of unstable modes nucleating at $A_c=0.889$.
However, we also find several smaller bubbles of instability at lower
values of $A_c$.  The first one we observe occurs near $A_c=0.333$,
which is shown magnified in Figs.~\ref{fig:mag:00} and
\ref{fig:arg:00}.  Together with Fig.~\ref{fig:kk:00}, these figures
show that this instability is generated by two nucleation events in
which the $\lkr=3$ and $\lkr=5$ eigenfrequencies collide, causing two
pairs of Floquet multipliers to leave the unit circle.  In more
detail, at $A_c=0$, there are four eigenfunctions of $E_T$ with
$\lkr=3$ and four with $\lkr=5$. The relevant eigenvalues are
\begin{equation*}
  \begin{aligned}
  \lambda_9 =
  \overline{\lambda_{10}} =
  \lambda_{11} =
  \overline{\lambda_{12}} =
  e^{-i\pi(2\sqrt{3}-4)} =
  e^{i\pi(0.535898)}, \\
  \lambda_{17} =
  \overline{\lambda_{18}} =
  \lambda_{19} =
  \overline{\lambda_{20}} =
  e^{i\pi(2\sqrt{5}-4)} =
  e^{i\pi(0.472136)}.
  \end{aligned}
\end{equation*}
We choose labels so the corresponding eigenfunctions $\dot q_j$ are
odd for $j\in\{9,10,17,18\}$ and even for $j\in\{11,12,19,20\}$.  As
$A_c$ increases over the range $0\le A_c\le0.31$, we see in
Figs.~\ref{fig:mag:00}, \ref{fig:arg:00} and~\ref{fig:kk:00} that the
eigenfrequencies $\sigma_9$ and $\sigma_{11}$ decrease, $\sigma_{17}$
and $\sigma_{19}$ increase, $\lkr_9$ and $\lkr_{11}$ remain nearly
equal to 3, and $\lkr_{17}$ and $\lkr_{19}$ remain nearly equal to 5.
Conjugate symmetry implies that $\sigma_{10}=-\sigma_9$,
$\lkr_{10}=\lkr_9$, etc., but the other pairs that are equal at
$A_c=0$ (e.g.~$\sigma_9=\sigma_{11}$ and $\lkr_9=\lkr_{11}$) split
apart slightly once $A_c>0$.  When $A_c$ reaches $0.3295$,
$\sigma_{9}$ collides with $\sigma_{17}$ (at $0.50813\pi$), and
$\lkr_{9}$ collides with $\lkr_{17}$ (at 4.1645).  This causes the
first pair of Floquet multipliers (with odd eigenfunctions) to leave
the unit circle. Shortly after this, at $A_c=0.3303$, $\sigma_{11}$
collides with $\sigma_{19}$ (at $0.50826\pi$), and $\lkr_{11}$
collides with $\lkr_{19}$ (at 4.1645).  This causes the second pair of
Floquet multipliers (with even eigenfunctions) to leave the unit
circle. The first pair of multipliers rejoins the unit circle at
$A_c=0.3362$, $\sigma_9=\sigma_{17}=0.50826\pi$,
$\lkr_9=\lkr_{17}=4.1569$, and the second pair rejoins at
$A_c=0.3372$, $\sigma_{11}=\sigma_{19}=0.50841\pi$,
$\lkr_{11}=\lkr_{19}=4.1568$, restoring stability to the standing
waves.  The mean wave numbers return to approximately 3 and 5 once
$A_c$ increases past 0.35 (see Fig.~\ref{fig:kk:00}).

The above scenario repeats itself for a number of different
eigenfrequency collisions between $0<A_c<0.889$.  However, the widths
of the intervals of $A_c$ over which the solutions are unstable are
quite small, so most standing waves in the range $0<A_c<0.889$ appear
to be stable. (Mercer and Roberts \cite{mercer:92} did not happen to
land on any of these windows of instability.)  In this parameter
regime, the mean wave number of each branch of Floquet multipliers
remains close to its initial integer value, except in regions of
instability, where the mean wave numbers of two branches briefly
coalesce while their eigenvalues lie outside the unit circle
(see Fig.~\ref{fig:kk:00}).

For larger values of crest acceleration ($A_c>0.889$), the situation
changes.  A number of nucleation events occur in which a Floquet
multiplier and its complex conjugate collide at $\lambda=1$ and split
off along the real axis. This leads to two new eigenvalues outside the
unit circle ($\{\lambda,\lambda^{-1}\}\subset\mathbb{R}$). By contrast,
for smaller values of $A_c$, four eigenvalues leave the unit circle
together ($\{\lambda, \overline\lambda, \lambda^{-1},
\overline{\lambda^{-1}}\} \subset\mathbb{C}$).  Rather than quickly
returning to the unit circle to regain stability, the magnitudes of many of
the unstable eigenvalues that nucleate with $A_c>0.889$ continue to grow as crest acceleration increases
(see Fig~\ref{fig:mag:00}). In addition, many more bubbles of
instability appear and disappear, with multiple nucleation events occurring in rapid
succession. In summary, all standing waves with $A_c>0.889$
appear to be unstable, with both the number of unstable modes and the
magnitude of the largest multiplier growing as $A_c$ increases.

\begin{figure*}
\includegraphics[width=.9\linewidth]{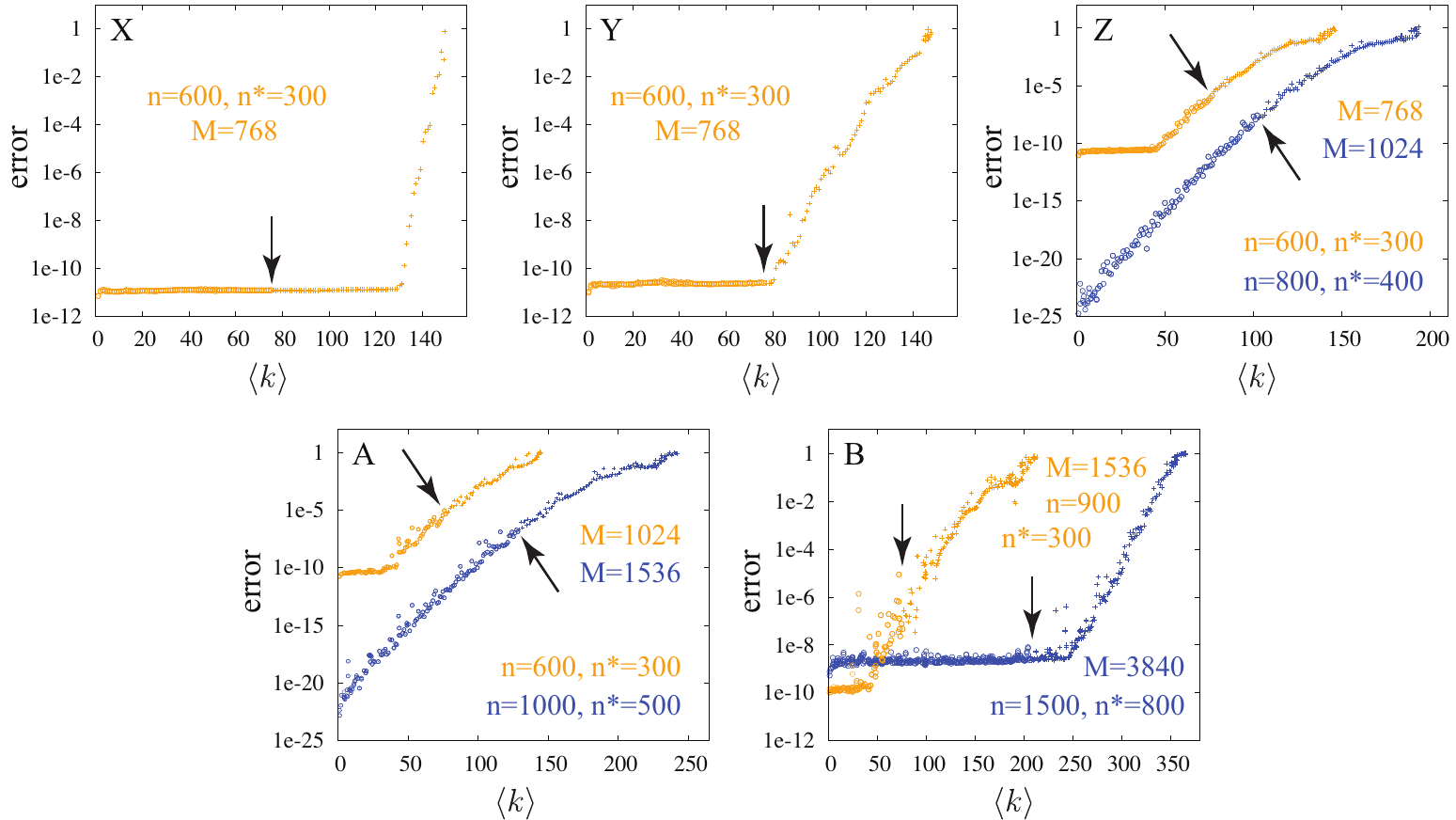}
\caption{\label{fig:eig:error}
  Plots of the residual error of an eigenvector versus its mean wave
  number for five standing wave solutions in deep water.  Arrows point
  to mode $n^*$, the last one used in the Floquet analysis of
  Section~\ref{sec:stable}. A second calculation was done on a finer
  grid for solutions Z, A, B to confirm accuracy. These were done in
  quadruple precision for Z and A, which is why the blue curves do not
  flatten out due to roundoff errors for small $\langle k\rangle$ in
  these two cases.  }
\end{figure*}

Next we check the accuracy of these Floquet multipliers.
Fig.~\ref{fig:eig:error} shows the residual error, defined in
(\ref{eq:err}), versus the mean wave number $\lkr$, for solutions X,
Y, Z, A, B above.  This residual error is made up of two parts:
roundoff error, which prevents $J_{1:n,1:n}z=\lambda z$ from holding
exactly, and truncation error, which arises when eigenvectors of
$J_{1:n,1:n}$ fail to belong to $\ker(J_{n+1:\infty,1:n})$
exactly. Roundoff error is further broken up into errors that occur
while solving the PDE to compute the matrix entries of $\hat E_T$, and
errors that occur in the numerical linear algebra routines while
diagonalizing $\hat E_T$. By posing the problem as an overdetermined
shooting method and monitoring the decay of Fourier modes
\cite{water2}, we can ensure that $\Delta t$ and $\Delta x$ are small
enough that truncation error in time-stepping is smaller than
roundoff errors.

The orange markers in Fig.~\ref{fig:eig:error} were computed with the
same parameters as the other solutions in
Figs.~\ref{fig:mag:00}--\ref{fig:arg:00}, given in the first three
rows of Table~\ref{tab:params}. The blue markers were computed with
the parameters in the remaining three rows of Table~\ref{tab:params}.
For small values of $\lkr$, roundoff error dominates truncation error
in these computations, leading to a flat region at the beginning of
each error plot. The two exceptions are the recomputation of solutions
Z and A in quadruple precision, where the time-periodic solution was
computed to roundoff-error accuracy \cite{water2}, but $M$ and $n$
were not quite large enough to reach roundoff error in the stability
study. We also see in these plots that increasing $M$ and $n$ causes
the truncation part of the residual error to shift downward, but can
increase the effect of roundoff error as more floating point operations
are involved.

\begin{figure}
\begin{center}
\includegraphics[width=\linewidth]{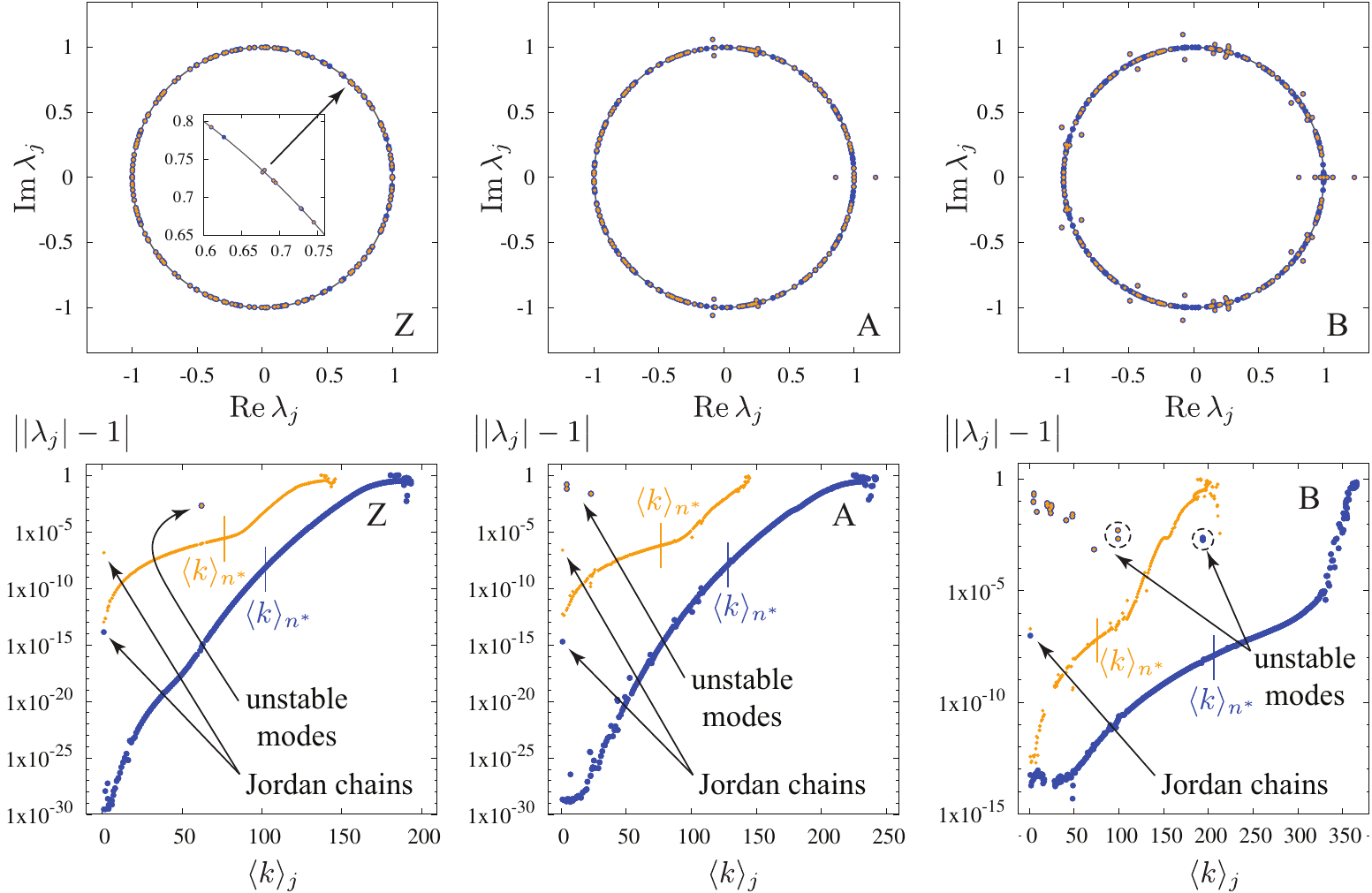}
\end{center}
\caption{\label{fig:stable:00} Plots of the first $n^*$ Floquet
  multipliers in the complex plane for solutions Z, A and B (top row)
  along with deviation of the first $n$ Floquet multipliers from the
  unit circle (bottom row).  This data is the same as in
  Fig.~\ref{fig:eig:error}, with blue markers corresponding to the
  second calculation on a finer grid. The dashed circles contain
  unstable eigenvalues that were missed on the coarse grid.}
\end{figure}

In the top row of Fig.~\ref{fig:stable:00}, we compare the first $n^*$
Floquet multipliers of solutions Z, A, B on the original mesh (orange
  markers, $n^*=300$) to the refined mesh (blue markers, $n^*=400$,
  $500$, $800$).  Each of the blue markers whose deviation from the
unit circle is visible at this resolution has an orange marker
directly on top of it, indicating that the original computations are
accurate enough to capture the most unstable modes. In the bottom row
of Fig.~\ref{fig:stable:00}, we plot the deviation of all $n$ computed
eigenvalues from the unit circle versus the mean wave number of the
eigenvalue. For the 6 solutions shown (Z coarse, Z, fine, A coarse, A
  fine, B coarse, B fine), we have $n=600, 800, 600, 1000, 900, 1500$,
respectively; see Table~\ref{tab:params} above.  A vertical line is
drawn on each curve at the cutoff $\langle k\rangle_j$ with $j=n^*$,
which separates the leading eigenvalues from the discarded
eigenvalues. For solution B, we see there are two additional clusters
of unstable modes (inside the dashed circles) that were not among the
first 300 Floquet multipliers. The first cluster was computed
accurately on the coarse mesh but was discarded due to $n^*$ being too
small, and the second cluster was not resolved on the coarse mesh but
emerged on the finer mesh.  These modes are stable enough that their
deviation from the unit circle in the upper right plot in
Fig.~\ref{fig:stable:00} is not visible. They can be seen in the
bottom right plot due to the logarithmic scale.

The markers labeled ``Jordan chains'' correspond to the eigenvalue
$\lambda=1$, which has algebraic multiplicity 4 but geometric
multiplicity 2 in the generic case. As explained in
Section~\ref{sec:jchains} below, there are two linearly independent
Jordan chains associated with $\lambda=1$, each of length 2.  This
causes errors of order $O(\epsilon^{1/2})$ when $\lambda=1$ splits
into two pairs of nearly reciprocal eigenvalues in floating-point
arithmetic, where $\epsilon=2.2\times10^{-16}$ in double precision and
$\epsilon=4.8\times10^{-32}$ in quadruple precision
\cite{Achain}. Each reciprocal pair of eigenvalues can split along the
real axis or onto the unit circle as a conjugate pair.  If either pair
splits along the real axis, the error appears as a spurious unstable
mode that deviates from the unit circle by roughly $10^{-7}$ in
double-precision or $10^{-14}$ in quadruple-precision. We note that
clusters of two or four Floquet multipliers of the form
$\{\lambda,\lambda^{-1}\}$ or
$\{\lambda,\lambda^{-1},\bar\lambda, \bar\lambda^{-1}\}$
appear as a single point in the bottom
plots of Fig.~\ref{fig:stable:00} since conjugation does not affect
the distance to the unit circle, and if $|\lambda|=1-\delta$ with
$\delta$ small, then $|\lambda^{-1}|=1+\delta+O(\delta^2)$; thus,
$\big||\lambda|-1\big|=|\delta|\approx\big||\lambda^{-1}|-1\big|$.

\subsection{Jordan chains associated with $\lambda=1$}
\label{sec:jchains}

The monodromy operator $E_T$ possesses at least two Floquet
multipliers equal to 1 since the water wave equations are invariant
under translation in time and space.  Indeed, if $q=(\eta,\varphi)$ is
a time-periodic solution of $\partial_t q = F(q)$, then
\begin{equation}\label{eq:p0p1}
  p_0(x,t)=\partial_tq(x,t) \qquad \text{and} \qquad
  p_1(x,t)=\partial_xq(x,t)
\end{equation}
satisfy the linearized equation $\partial_t p=DF(q(t))p$ as well as
$p(x,T)=p(x,0)$.  Thus,
\begin{equation}
  E_T [p_j(\cdot,0)] = \lambda p_j(\cdot,0), \qquad \lambda=1, \;\; j=0,1.
\end{equation}
For symmetric standing waves, $\eta$ and $\varphi$ are even functions
of $x$ for all $t$, so $p_0(x,0)$ is even while $p_1(x,0)$ is odd.

There are three ways to compute the eigenfunctions $p_0(x,0)$ and
$p_1(x,0)$.  We can differentiate the nonlinear solution $q$ with
respect to $t$ and $x$, respectively, as in (\ref{eq:p0p1}).  We can compute all
the eigenvalues and eigenfunctions of a numerical approximation of
$E_T$, splitting them into even and odd eigenfunctions, as was done in
Figures~\ref{fig:mag:00}--\ref{fig:stable:00} above. Extracting the
$\lambda=1$ eigenfunctions leads to two copies of $p_0(x,0)$ and two copies of $p_1(x,0)$,
up to constant factors. (Jordan chains are responsible
  for the duplicate eigenvectors, which cause the eigenvector
  matrix to be singular.) Or we can
compute the kernel of $E_T-I$.  This last approach has the advantage
that with minor modification, we can obtain Jordan chain information
as well.  Using the algorithm of Wilkening \cite{Achain}, we have
found that $p_0(x,0)$ and $p_1(x,0)$ generally have an associated
vector $p_j^\e1(x,0)$ such that
\begin{equation}\label{eq:jordan1}
  E_T [p_j^\e1(\cdot,0)] = p_j(\cdot,0) + \lambda p_j^\e1(\cdot,0), \qquad \lambda=1.
\end{equation}
These associated vectors have a natural physical interpretation.  For
example, if we think of the period as a bifurcation parameter
controlling $q(x,t;T)$, then time-periodicity gives
$q(x,T;T)=q(x,0;T)$.  Thus, $\partial_tq(x,T;T) + \partial_Tq(x,T;T)
= \partial_Tq(x,0;T)$.  It follows that
\begin{equation}\label{eq:p01:def}
  p_0^\e1(x,t) = -\partial_Tq(x,t;T)
\end{equation}
is a solution of the linearized equation and satisfies
\begin{equation}\label{eq:ET:jchain}
  E_T[p_0^\e1(\cdot,0)] = p_0^\e1(\cdot,T) = p_0(\cdot,0) + p_0^\e1(\cdot,0),
\end{equation}
as required. Iterating this equation gives
\begin{equation}\label{eq:jchain:m}
  E_T^m[p_0^\e1(\cdot,0)] = mp_0(\cdot,0) + p_0^\e1(\cdot,0), \qquad (m\in\mbb Z),
\end{equation}
which can be interpreted as expressing
\begin{equation}\label{eq:jchain:eps}
  \begin{aligned}
    q(x,mT;T-\veps) &\approx q(x,m(T-\veps);T-\veps) + \veps m q_t(x,m(T-\veps);T-\veps) \\
    & = q(x,0;T-\veps) + \veps m q_t(x,0;T-\veps) \\
    & \approx q(x,0;T) + \veps[m q_t(x,0;T) - q_T(x,0;T)], \qquad (m\in\mbb Z).
  \end{aligned}
\end{equation}
The idea here is that if the initial condition is perturbed in the
direction of another standing wave with a slightly shorter period,
then each time $t$ advances by the period $T$ of the original
solution, it will travel slightly further in the $q_t(\cdot,0)$
direction. We note that these Jordan chains lead to ``secular
instabilities'' in the solution that grow linearly in time. We will
still consider a solution to be stable if these are the only
instabilities.

This argument that Jordan chains of the form (\ref{eq:p01:def}) should
exist assumes that standing waves occur in smooth families.  Even if
this is not the case, e.g.~if there are infinitely many disconnections
of the type shown in Figure~\ref{fig:bif:00:jump} on scales too small
to be observed in quadruple-precision arithmetic, the monodromy
operator may still have Jordan chains. In our numerical experiments,
these Jordan chains are always present, except at turning points,
where they become genuine eigenvectors. For example, $T$ reaches a
relative maximum of $T=6.543\,698$ at $A_c=0.889\,283$.  Near such a
turning point, it is better to use a different coefficient $c_k$ in
(\ref{eq:init:stand}) than $c_0=T$ as the bifurcation parameter; we
used $c_5$.  From $q(x,T;c_k)=q(x,0;c_k)$ we have
\begin{equation}
  \der{T}{c_k}\der{q}{t}(x,T;c_k) + \der{q}{c_k}(x,T;c_k) = \der{q}{c_k}(x,0;c_k).
\end{equation}
Thus, rescaling the associated vector via $p_0^\e1(x,t)=-(\partial
  q/\partial c_k)(x,t;c_k)$, the first term on the right-hand side of
(\ref{eq:jordan1}) acquires a factor of $\pa T/\pa c_k$ that vanishes
at a turning point. As a result, $p_0^\e1(x,t)$ becomes an independent eigenvector
rather than blowing up like $\partial q/\partial T$ does. Similarly,
using $c_1$ to parametrize the bifurcation from the zero-amplitude wave, we have
$\pa T/\pa{c_1}=0$, which explains why Jordan chains did not arise in
Section~\ref{sec:zro:stab}.

We expected to find the associated vector $p_0^\e1(\cdot,0)$ when we
computed Jordan chains at $\lambda=1$ using the numerical algorithm of
\cite{Achain}, but we did not anticipate finding the second associated
vector $p_1^\e1(\cdot,0)$. Just as perturbation in the direction of
$p_0^\e1(\cdot,0)$ leads to linear drift in time over multiple cycles
due to a slight change in period, perturbation in the
$p_1^\e1(\cdot,0)$ direction should lead to linear drift in space.
This suggests the existence of traveling-standing waves bifurcating
from the pure standing wave solutions. We have investigated this, and
have successfully computed a two-parameter family of
traveling-standing waves, which will be reported on elsewhere
\cite{waterTS}.

\subsection{Counter-propagating solitary waves in shallow water}
\label{sec:shallow}

\begin{figure}
\begin{center}
\includegraphics[width=.92\linewidth]{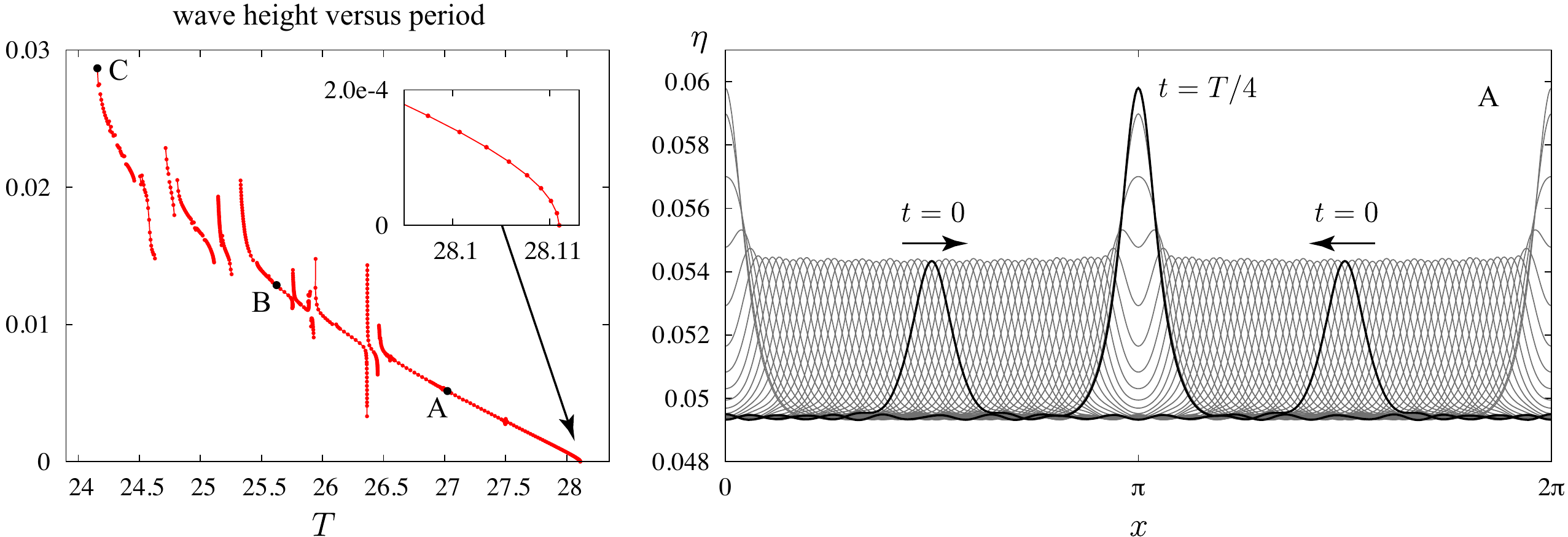}
\end{center}
\caption{\label{fig:bif:005} Wave height versus period for a family of
  standing waves in shallow water ($h=0.05$) and time-elapsed
  snapshots of solution A over a period. (left) Many more
  disconnections are visible at this depth than were observed in the
  infinite-depth case above. (right) The arrows show the initial
  direction of motion at $t=0$. The waves return to
  this initial position at $t=T/2$, but moving in the opposite
  direction.  Small vertical oscillations are visible in the wave
  troughs and crest trajectories as the waves travel back and forth
  between collisions.  }
\end{figure}

Wilkening and Yu \cite{water2} found that varying the fluid depth of
standing water waves leads to nucleation or vanishing of loop
structures in the bifurcation curves that meet or nearly meet at
perfect or imperfect bifurcations.  Closed loops can even nucleate in
isolation and later join with pre-existing families of
solutions. Disconnections and bifurcations tend to disappear as the
fluid depth increases, though some persist in the infinte depth limit,
as seen above in Fig.~\ref{fig:bif:00:jump}.  These disconnections and
gaps in the bifurcation curves can be interpreted \cite{water2} as
numerical manifestations of the Cantor-like structures that arise in
analytical studies of standing water waves due to small divisors
\cite{plotnikov01, iooss05}.  In shallow water, resonances abound, and
the bifurcation curves are highly fragmented even in
double-precision. This is illustrated in the left panel of
Fig.~\ref{fig:bif:005}, which shows wave height versus period for
symmetric standing waves of wavelength $2\pi$ in a fluid of depth
$h=0.05$. Here wave height is defined as the crest to trough height at
$t=T/4$, when the fluid comes to rest. This differs by a factor of two
from the convention used in Section~\ref{sec:deep} above.  These
solutions were first computed in \cite{water2}, where a different set
of bifurcation curves are given that show how $\hat\varphi_k(0)$
varies with $T$ for $k=1$, 17 and 36. The periods of these waves are
much longer than their infinite depth counterparts. For example, in
the linear regime, $T=2\pi/\sqrt{gK \tanh K h}=28.1$ when $h=0.05$
and $T=6.28$ when $h=\infty$. In both cases, $g=1$ and $K=1$.

\begin{figure}
\begin{center}
\includegraphics[width=.67\linewidth]{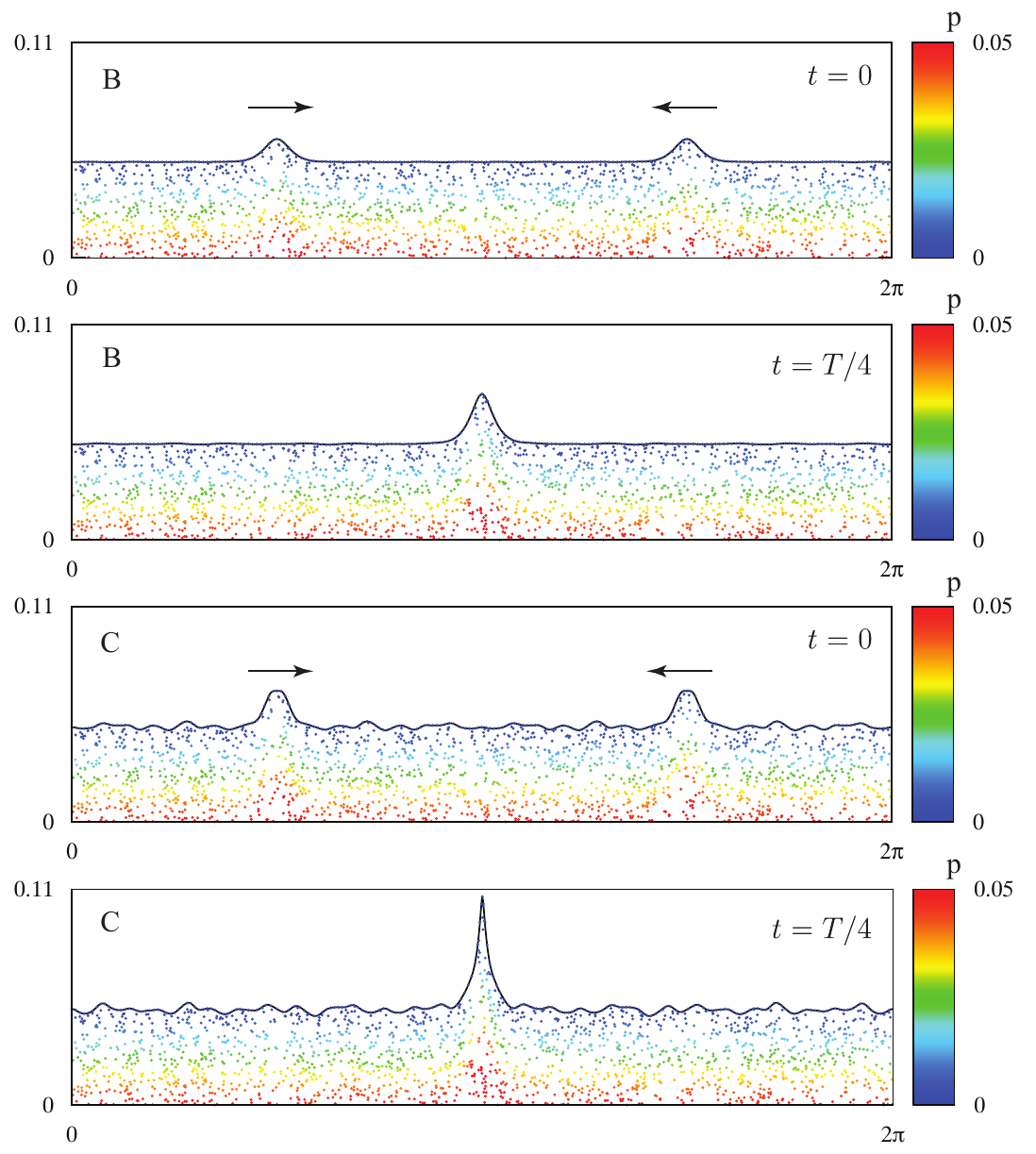}
\end{center}
\caption{\label{fig:collision:BC} Solutions B
  and C from Figure~\ref{fig:bif:005} at times $t=0$ and $t=T/4$.
  Solution B remains calm between collisions while solution
  C is choppy.  These images are taken from movies in which passively
  advected particles have been added to the fluid for visualization,
  color coded by pressure. }
\end{figure}

Beyond the linear regime, standing waves in shallow water take the
form of counter-propagating solitary waves that repeatedly collide
with one another. These waves are special in the sense that
low-amplitude, high-frequency oscillatory modes are present initially
and have been tuned so that no new radiation is generated each time
the solitary waves collide. The right panel of Fig.~\ref{fig:bif:005}
shows time-elapsed snapshots of solution A in the bifurcation diagram
in the left panel, and Fig.~\ref{fig:collision:BC} shows solutions B
and C at $t=0$ and $t=T/4$. Figure~\ref{fig:collision:BC} also shows
the pressure beneath the waves and marker particles for visualizing
the flow in movies of these simulations. Aside from the
solitary waves, solution B remains calm throughout its evolution.
By contrast, the solitary waves of solution C travel over a choppy
and erratic background of higher-frequency waves that are clearly
visible in static images.

\begin{figure}
\begin{center}
\includegraphics[width=.91\linewidth]{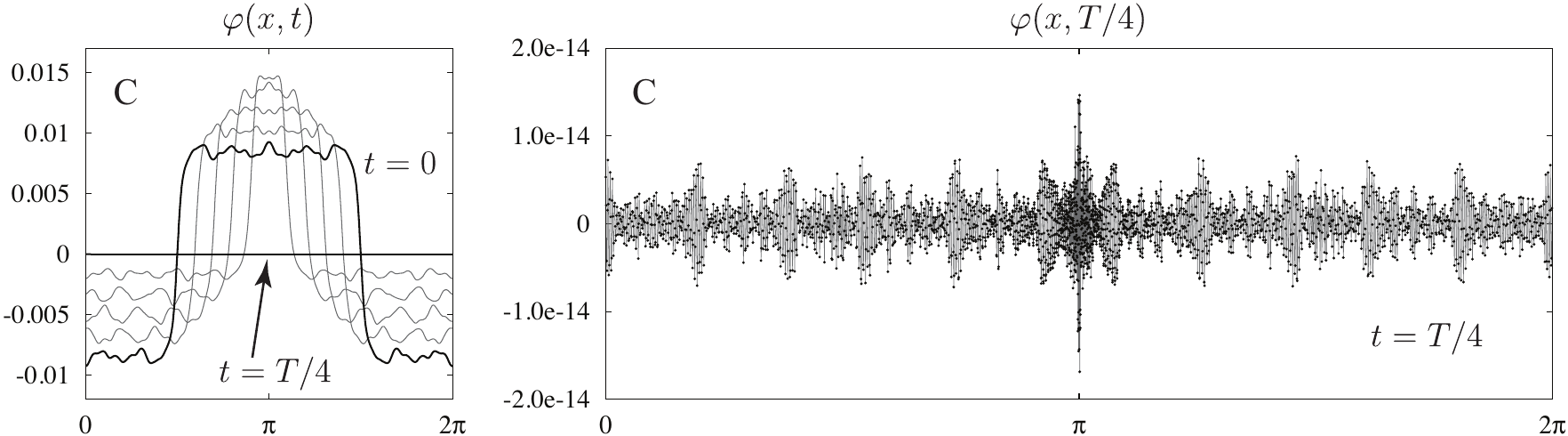}
\end{center}
\caption{\label{fig:error:005} Evolution of the velocity potential of
  solution C over a quarter period. (left) The velocity potential
  jumps rapidly at the location of the solitary waves, where fluid
  velocity $\mathbf{u}=\nabla\phi$ is highest.  (right) At $t=T/4$,
  the velocity potential has been driven to zero (up to roundoff
    error), yielding a minimum value of $f=3.3\times 10^{-30}$ for the
  objective function in (\ref{eq:f:phi}). For mesh refinement, we set
  $\nu=2$, $\theta_1=0.85$, $\theta_2=0.15$, $\gridkap_1=0$,
  $\gridkap_2=2$, $\rho_2=0.15$, $M_1=2048$, $N_1=204$, $M_2=3072$ and
  $N_2=108$.  The computation took 43 minutes using an NVIDIA Tesla
  M2050 GPU, with 32 minutes devoted to computing the Jacobian.}
\end{figure}

The height of the solitary waves relative to the fluid depth is
remarkable for both of these solutions.  For solution C, when the
solitary waves collide and the fluid comes to rest, the wave crest
reaches a height of $0.10636$ (measured from the bottom), which is
more than twice the fluid depth $h=0.05$. It is surprising that such a
jet could form and then return back along its path to form a
time-periodic solution. Figure~\ref{fig:error:005} shows snapshots of
$\varphi(x,t)$ for solution C and a plot of the residual error at
$t=T/4$, confirming that the wave actually comes to rest at this
time. Because the wave amplitudes of these solutions are so large, it
is not appropriate to consider $\eta(x,t)$ as a small perturbation of
a flat wave profile, to ignore vertical components of velocity, nor to
assume hydrostatic pressure forces. Thus, most of the assumptions one
makes in deriving model equations such as Serre's equations, the
Boussinesq system, or the Korteweg de-Vries equation, break down
\cite{vanden:broeck:book}. These wave interactions are fully (as
  opposed to weakly) nonlinear.  Nevertheless, solutions such as B
have many properties in common with the perfectly elastic collisions
of solitons that occur in integrable equations.

\begin{figure}
\begin{center}
\includegraphics[width=.95\linewidth]{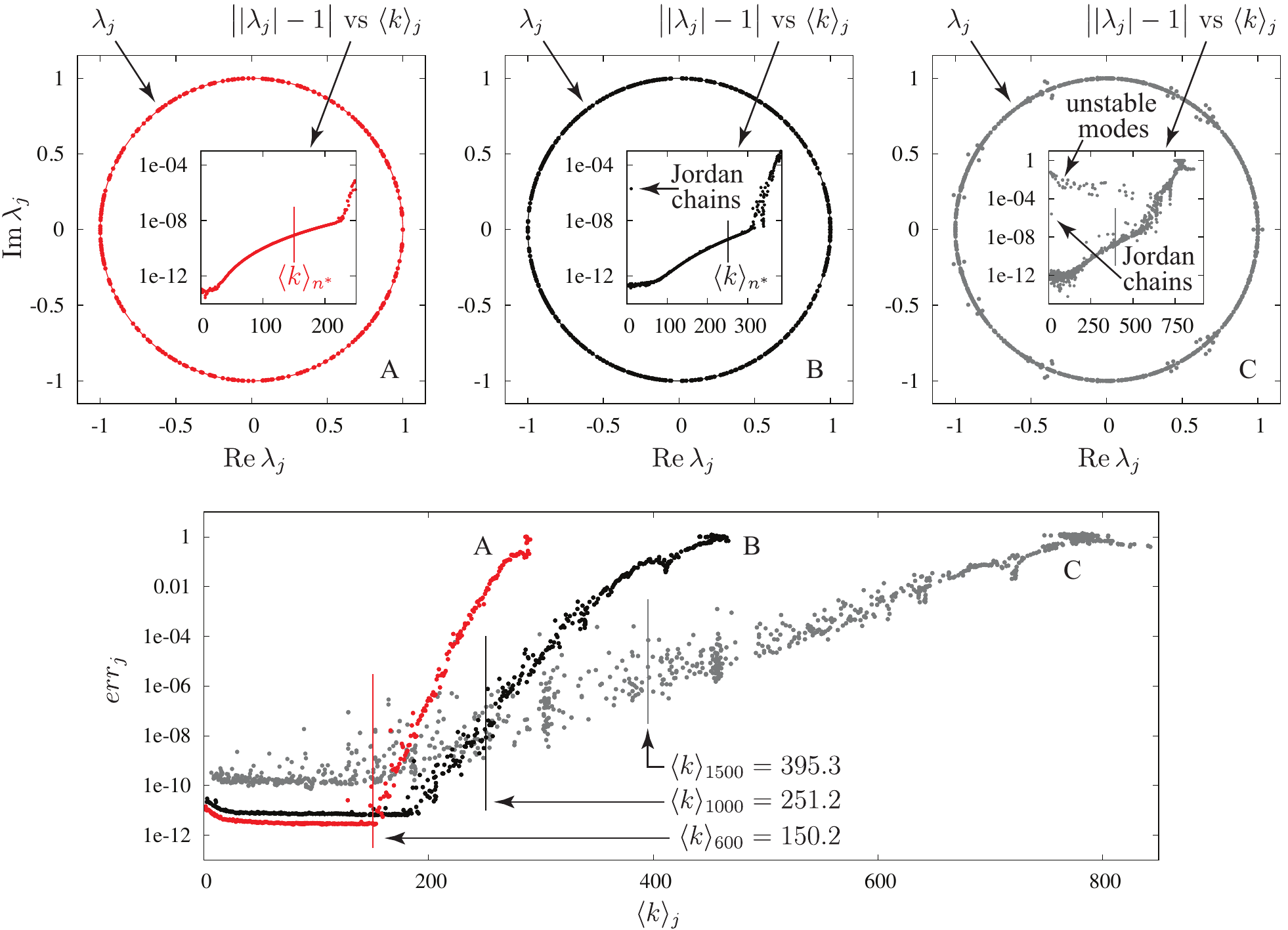}
\end{center}
\caption{\label{fig:stable:005} Floquet stability analysis of
  solutions A--C in Figures~\ref{fig:bif:005}--\ref{fig:error:005}.
  (top row) The first $n^*=600$, $1000$ and $1500$ Floquet multipliers
  $\lambda_j$ are plotted in
  the complex plane. The inset graphs show that the more reliable
  eigenvalues (of lower mean wave number) satisfy $|\lambda_j|=1$ with
  high accuracy, except for the unstable modes of solution C and the
  eigenvalue $\lambda=1$, which has Jordan chains that cause it to
  split. (bottom) Residual error versus mean wave number.  The vertical
  lines give the cutoff $\langle k\rangle_{n^*}$ beyond which the
  eigenvalues were discarded.  }
\end{figure}

Next we consider the stability of these solutions with respect to
harmonic perturbations. The results are summarized in
Figure~\ref{fig:stable:005}. For solutions A and B, we used uniform
spatial grids with $M^\e{A}=1536$ and $M^\e{B}=2048$ gridpoints.  For
solution C, we used mesh refinement near $(x,t)=(\pi,T/4)$ and
$(x,t)=(0,3T/4)$. More specifically, in the notation of
Section~\ref{sec:eqm}, we set $\nu=5$, $\gridkap_l=\{0,2,0,-2,0\}$,
$\theta_l=\{0.21,0.08,0.42,0.08,0.21\}$, and $M^\e{C}_l=\{3888$,
$5184$, $3888$, $5184$, $3888\}$.  For solutions $\{A,B,C\}$, we
computed $n=\{1200$, $2000$, $4000\}$ Floquet multipliers, kept the
first $n^*=\{600,1000,1500\}$ ordered by mean wave number, and
discarded the rest.  All the multipliers of solutions A and B lie on
the unit circle, with $\big||\lambda_j|-1\big|$ ranging from
$10^{-13}$ for low values of $\langle k\rangle$ to $10^{-9}$ for large
values of $\langle k\rangle$, up to the cutoff point $\langle
k\rangle_{n^*}$.  As in the infinite depth case, a pair of Jordan
chains of length 2 is associated with the two-dimensional eigenspace
at $\lambda=1$.  Numerically, this degenerate eigenvalue splits into a
cluster of four simple eigenvalues with errors on the order of
$\sqrt\epsilon$, where $\epsilon$ is machine precision
\cite{demmel}. These four eigenvalues consist of two pairs of nearly
reciprocal eigenvalues, i.e.~the product of the eigenvalues in each
pair differs from 1 by $O(\epsilon)$.  For solution A, the degenerate
eigenvalue happens to split into two complex conjugate pairs that lie
on the unit circle (since $\lambda^{-1}=\bar\lambda+O(\epsilon)$);
thus, they do not appear anomalous in the plot of
$\big||\lambda_j|-1\big|$.  For solution B, the degenerate eigenvalue
split into two pairs of real eigenvalues; and for solution C, the
degenerate eigenvalue split into one pair of real eigenvalues and one
pair of complex conjugate eigenvalues.  Regardless of how they split,
which is unpredictable due to roundoff errors, these Jordan chains
correspond to nearby families of time-periodic solutions and
traveling-standing waves. Thus, solutions A and B appear to be
linearly stable to harmonic perturbation (aside from the secular
  instabilities of the Jordan chains). In numerical experiments that will be
reported elsewhere \cite{water:LT}, we find that nearby initial
conditions evolved under the nonlinear equations of motion remain
close to these stable time-periodic solutions over thousands of
cycles.

By contrast, solution C has many eigenvalues $\lambda_j$ outside the
unit circle. Perturbation of the time-periodic solution in these
directions will lead to exponential growth.  The growth is actually
quite slow, with the largest multiplier being
$\rho=|\lambda_\text{max}|=1.062$. This is the amplification factor
over a cycle (consisting of two collisions of the solitary waves), not
the growth rate per unit time.  Thus, if roundoff error contributes an
error of $10^{-13}$ per cycle, the accumulated error after $p$ cycles
is expected to be around $10^{-13}\rho(\rho^p-1)/(\rho-1)$.  This is
roughly what happens in a long-time simulation, except there is a
large startup phase of about 100 cycles where the solution remains
time-periodic to 13 digits. The next 380 cycles look time-periodic to
the eye, but exhibit exponential growth when
$\|\eta(\cdot,T)-\eta(\cdot,0)\|$ is computed.  Over the next 90
cycles, the solitary waves fall out of phase with each other and the
background radiation grows in amplitude.  Finally, after evolving
solution C through 571 cycles, the numerical solution blows up, with
GMRES failing to converge when solving for the dipole density $\mu$.

\ignore{
\begin{figure}[b]
\begin{center}
\includegraphics[width=\linewidth]{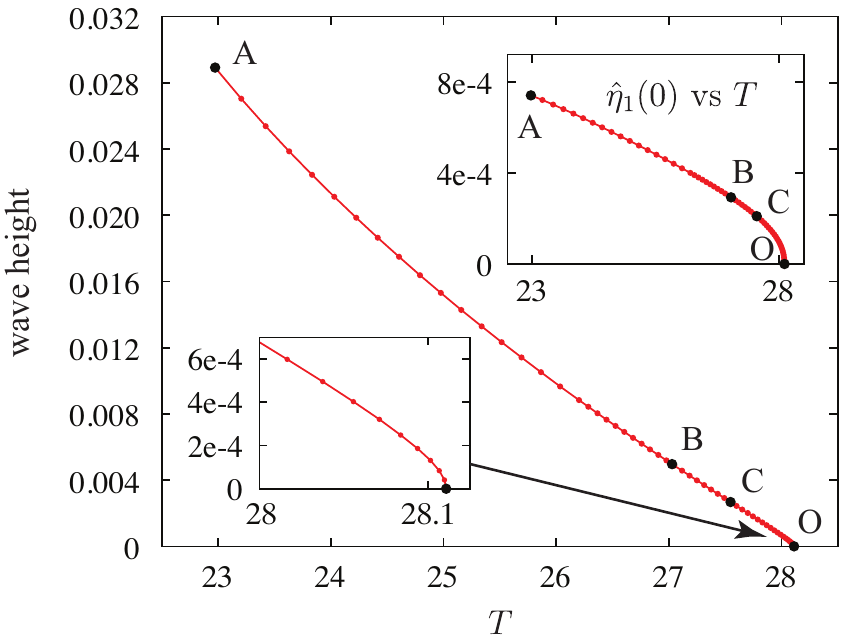}
\end{center}
\caption{\label{fig:bif:stokes} Collision of two right-moving Stokes
  waves that nearly return to their initial configuration after the
  interaction. (left) Plots of wave height and $\hat\eta_1(0)$ versus
  period for spatially $2\pi$-periodic Stokes waves.  The temporal
  periods of solutions A and C are $6T_A=137.843\approx 137.738 =
  5T_C$.  (center) Solutions A and C were combined via
  (\ref{eq:AandC}) and evolved through one collision to $t=137.738$.
  (right) Through trial and error, we adjusted the amplitude of the
  smaller Stokes wave and the simulation time to obtain a nearly
  time-periodic solution.  }
\end{figure}
}

\subsection{Gravity-capillary solitary wave interactions in deep water}
\label{sec:surf}

\begin{figure}
\begin{center}
\includegraphics[width=\linewidth]{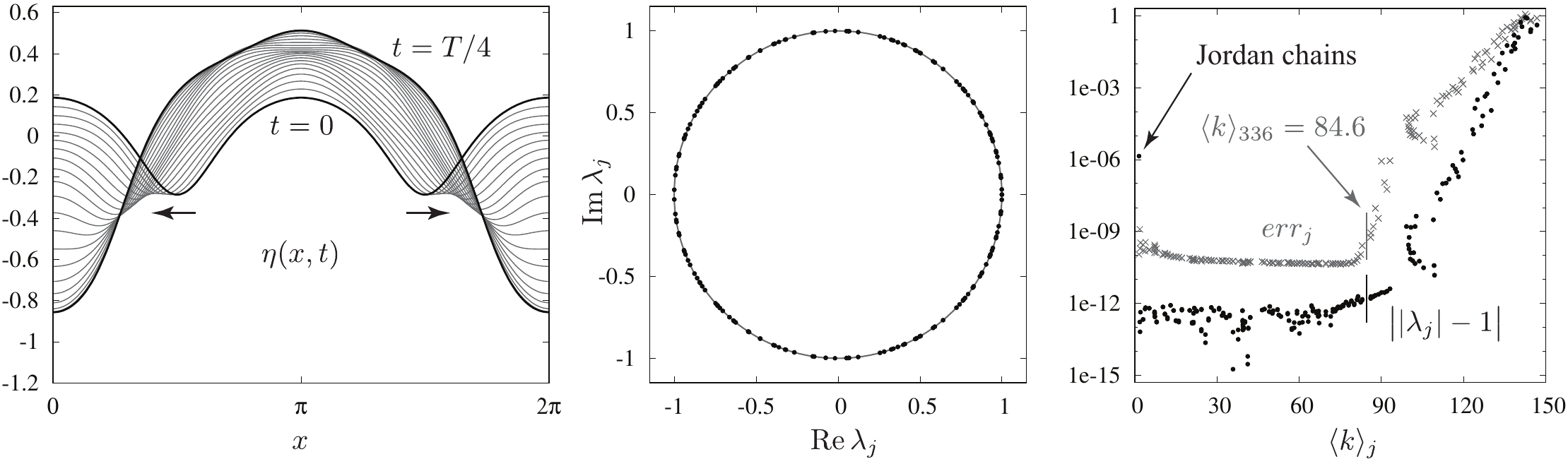}
\end{center}
\caption{\label{fig:concus} Stability calculation of a
  moderate-amplitude standing gravity-capillary water wave with
  crest-to-trough height of 1.3693. (left) At this amplitude, the wave
  consists of counter-propagating depression waves that repeatedly
  collide and change direction. (center) The leading $n^*=336$ Floquet
  multipliers lie on the unit circle. (right) Plots of $err_j$ and
  $\big||\lambda_j|-1\big|$ versus mean wave number.  We computed
  $n=600$ eigenvalues of the monodromy operator on a grid with
  $M=1024$ gridpoints and $N=28800$ timesteps over a full period. The
  underlying wave only needs $M=256$ gridpoints and 720 timesteps over
  a quarter period to drive the objective function $f$ in
  (\ref{eq:f:phi}) down to $3.7\times10^{-28}$, but additional
  timesteps are needed to resolve the dynamics of high-frequency
  perturbations. Surface tension increases the stiffness of the water
  wave equations, and requires more timesteps than the $\sigma=0$ case
  with an explicit method. We did not employ the small scale
  decomposition strategy of Hou, Lowengrub and Shelley \cite{hls94}.}
\end{figure}

We now investigate the effect of surface tension on the dynamics and
stability of time-periodic water waves.  For simplicity, we consider
only the infinite depth case. Concus \cite{concus:62} and
Vanden-Broeck \cite{vandenBroeck:84} computed the leading terms in a
perturbation expansion for standing water waves with surface tension,
building on the work of Penney and Price \cite{penney:52} and
Tadjbakhsh and Keller \cite{tadjbakhsh}.  Concus predicted that the
quadratic correction to the period in the infinte depth case is
\begin{equation}\label{eq:concus}
  T = \frac{2\pi}{\sqrt{1+\gamma}}\left(
    1+\frac{8-27\delta-36\delta^2-81\delta^3}{32(1-9\delta^2)}\epsilon^2 +
    O(\epsilon^4)\right), \qquad \epsilon = c_1 = \hat\varphi_1(0),
\end{equation}
where $\gamma=\frac{\sigma k^2}{\rho g}$ is a dimensionless surface
tension parameter, $\delta=\frac{\gamma}{1+\gamma}$ is the
relative capillarity, and $c_1$ refers to the notation of Equation
(\ref{eq:init:stand}).  Wilkening and Yu \cite{water2} computed a
family of standing waves of this type using the overdetermined
shooting method described in Section~\ref{sec:osm} and confirmed
(\ref{eq:concus}) in the special case of $\gamma=1$ and
$\delta=1/2$. It was also found in \cite{water2} that beyond the
linear regime, larger-amplitude standing waves of this type take the
form of counter-propagating solitary depression waves that repeatedly
collide and reverse direction.  We computed the leading Floquet
multipliers of a moderate-amplitude wave in this family, namely
solution B of Section 4.6 of \cite{water2}, and found that it is
stable to harmonic perturbations. The results are summarized in
Figure~\ref{fig:concus}.

Next we search for a new type of time-periodic gravity-capillary water
wave built from counter-propagating solitary waves of the type
discovered by Longuet-Higgins \cite{lh89} and studied by Vanden Broeck
\cite{vandenBroeck:cap} and Milewski et.~al.~\cite{milewski:11}. Our
idea is to collide two identical traveling waves of this type (moving in
  opposite directions) together, see how they interact, and optimize
the initial conditions of the combined wave to obtain a time-periodic
solution.

The first step is to compute traveling waves. While several methods
have been developed previously for this purpose \cite{rienecker,
  schwartz:fenton,vandenBroeck:cap,milewski:11,wilton:16}, we modified the
trust-region code we wrote to compute standing water waves so that it
can also compute traveling waves. Details are given in
Appendix~\ref{sec:trav}.  Instead of computing entire families of
traveling gravity-capillary waves, we will focus on two solutions,
with amplitudes close to the ones plotted in Figure~4 of
\cite{milewski:11}. Following the conventions of that paper, we
temporarily assume $\sigma/\rho=1$ and $g=1$ in (\ref{eq:ww}) and try
starting guesses of the form
\begin{equation}\label{eq:trav:cap:guess}
  \begin{aligned}
    &\text{Solution A:} \qquad \eta_0(x) = -0.15\frac{\cos x - a}{1-a}
    \opn{sech}(x/6.4), \qquad \varphi_0 = c H\eta_0, \qquad
    \left\{
    \begin{aligned}
      a &= \opn{sech}(3.2\pi), \\
      c &= 1.408,
      \end{aligned}\right. \\
    &\text{Solution B:} \qquad \eta_0(x) = -0.32\frac{\cos x - a}{1-a}
    \opn{sech}(x/2.4), \qquad \varphi_0 = c H\eta_0, \qquad \left\{
    \begin{aligned}
      a &= \opn{sech}(1.2\pi), \\
      c &= 1.385.
    \end{aligned}\right.
  \end{aligned}
\end{equation}
These functions $\eta_0$ have roughly the right oscillation
frequencies and decay rates (within 20 percent) as the graphs given in
\cite{milewski:11}.  We have estimated $\varphi$ from the linearized
equations (\ref{eq:lin}), which imply for a wave traveling at speed
$c$ that $-c\eta_x = \mc{G}\varphi$, or $\varphi = cH\eta$, where $H$
is the Hilbert transform.  While linear theory is not an accurate
model in this regime, it works well enough as a starting point for the
trust region shooting method to converge to a traveling solution of
the nonlinear equations.  The values of wave speed, $c$, were
estimated from the graph in Figure~1 of \cite{milewski:11}, and the
parameter $a$ was chosen so that $\eta_0(x)$ has zero mean.  The
actual values of $c$ for waves of amplitude $-0.15$ and $-0.32$
turned out to be $1.40995$ and $1.38558$, respectively.

In our code, we assume the domain is $[0,2\pi]$; however, the
convention that $\sigma/\rho=1$ leads to traveling waves that extend
well outside of this range.  So instead of choosing
$L=(\tilde\sigma/\tilde\rho \tilde g)^{1/2}$ and $\tau=(\tilde
  \sigma/\tilde \rho \tilde g^3)^{1/4}$ as characteristic length and
time scales, with notation as in (\ref{eq:nondim}), we use
$L=40(\tilde\sigma/\tilde\rho \tilde g)^{1/2}$ and
$\tau=\sqrt{40}(\tilde\sigma/\tilde\rho \tilde g^3)^{1/4}$. For real
water (assuming $\tilde\sigma=72 \text{ dyne}/\text{cm}$), this is
$L=10.8$ cm and $\tau=0.105$~s.  In the dimensionless equation
(\ref{eq:ww}), this causes $\sigma/\rho$ to change from 1 to
$0.000625$ and $g$ to remain~1.  The starting guesses
(\ref{eq:trav:cap:guess}) are modified by replacing $x$ by $40x$ on
the right-hand-side in the formulas for $\eta_0(x)$, dividing $0.15$
and $0.32$ by 40, and dividing $1.408$ and $1.385$ by $\sqrt{40}$.
The numerical values of wave height and velocity potential will then
decrease by factors of $40^{-1}$ and $40^{-3/2}$, respectively.  With
starting guess A for the traveling wave, the trust region shooting
method reduced $f$ to $3.8\times 10^{-24}$ in only 12 function
evaluations and 2 Jacobian calculations (using $M=2048$ gridpoints and
  $n=901$ unknowns).  With starting guess B, 10 function evaluations
and 1 Jacobian calculation were sufficient to minimize $f$ to
$4.5\times 10^{-26}$, using $M=3072$ gridpoints and $n=1801$ unknowns.
The former calculation took 34 seconds while the latter took 111
seconds.  In both cases, the amplitude, $\eta(0,0)$, was fixed via
(\ref{eq:eta:constr}) to be $-0.15/40$ for solution A and $-0.32/40$
for solution B.

Once a traveling wave is found, we obtain a starting guess for
counter-propagating time-periodic solutions by defining the initial
conditions
\begin{equation}\label{eq:counter:prop:trav}
  \begin{aligned}
  \eta_0(x) &= \eta_0^\e{\text{trav}}(x-\pi/2) +
              \eta_0^\e{\text{trav}}(x-3\pi/2), \\
  \varphi_0(x) &= \varphi_0^\e{\text{trav}}(x-\pi/2) -
                 \varphi_0^\e{\text{trav}}(x-3\pi/2).
  \end{aligned}
\end{equation}
The sign change in velocity potential causes the second wave to
travel left. The Fourier modes of the initial conditions
take the form (\ref{eq:init:trav}) for traveling waves and
(\ref{eq:init:stand}) for counter-propagating waves.  They are
related via
\begin{equation}
  \hat\eta_k(0) = 2\cos(k\pi/2)\,\hat\eta^\e{\text{trav}}_k(0), \qquad
  \hat\varphi_k(0) = -2i\sin(k\pi/2)\,\hat\eta^\e{\text{trav}}_k(0), \qquad
  (k\in\mathbb{Z}).
\end{equation}
From this starting guess, we minimized the objective function
(\ref{eq:f:phi}) assuming initial conditions of the form
(\ref{eq:init:stand}). We specified $\eta(\pi/2,0)$ to be $-0.00375$
for solution A and $-0.008$ for solution B by appending equation
(\ref{eq:eta:constr}) to the residual vector $r$.

\begin{figure}
\begin{center}
\includegraphics[width=\linewidth]{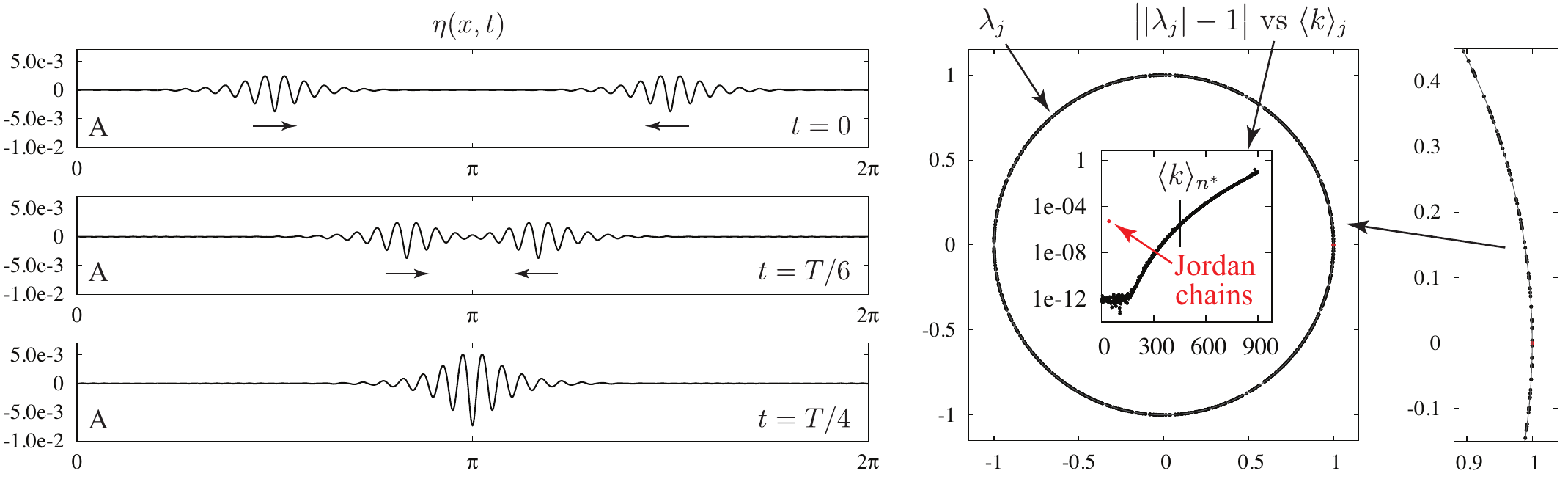}
\end{center}
\caption{\label{fig:cap:705} Counter-propagating gravity-capillary
  waves with $\sigma/\rho=0.000625$ and initial amplitude
  $\eta(\pi/2,0)=\eta(3\pi/2,0)=-0.00375$. At time $t=T/4$, the
  velocity potential is zero and the fluid comes to rest.  The period
  is $T=28.206$. We used $M=2048$ grid points and $n=620$ initial
  unknown Fourier modes to minimize $f$ to $2.41\times 10^{-30}$ with
  only 8 function evaluations and 3 Jacobian calculations. (right) The
  first $n^*=1800$ Floquet multipliers all lie on the unit circle,
  indicating that this solution is linearly stable to harmonic
  perturbations. We computed $n=3600$ columns of $\hat E_T$ with
  $M=3600$ gridpoints and $N=2700$ timesteps. The calculation took 5
  days. As in Section~\ref{sec:shallow}, the eigenvalue $\lambda=1$
  has two linearly independent Jordan chains of length two that cause
  the eigenvalue to split in the numerical calculation of Floquet
  multipliers. }
\end{figure}

\begin{figure}
\begin{center}
\includegraphics[width=\linewidth]{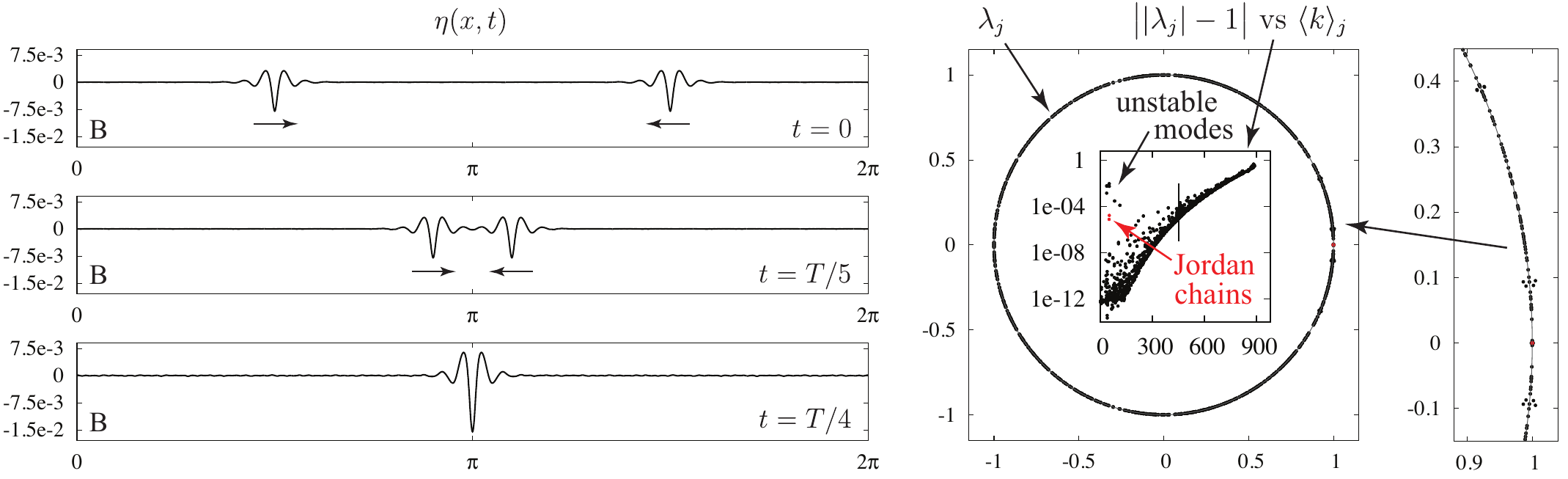}
\end{center}
\caption{\label{fig:cap:525} Counter-propagating gravity-capillary
  waves with $\sigma/\rho=0.000625$ and initial amplitude
  $\eta(\pi/2,0)=\eta(3\pi/2,0)=-0.008$. We used $M=3072$ gridpoints
  and $n=900$ unknown Fourier modes of the initial condition to obtain
  $f=3.9\times 10^{-28}$. The period of this solution is $T=28.711$.
  (right) Some of the first 1800 Floquet multipliers lie outside the
  unit circle, indicating that this solution is linearly unstable.
}
\end{figure}

The results are summarized in Figures~\ref{fig:cap:705}
and~\ref{fig:cap:525}.  In both cases, the wave packets approach each
other without changing shape until they start to overlap.  As they
collide, they produce a localized standing wave that grows in
amplitude to the point that $\varphi$ becomes exactly 0 (at time
  $T/4$).  The standing wave then decreases in amplitude and the two
traveling waves emerge and depart from one another.  As in
Section~\ref{sec:shallow} above, the background radiation is
synchronized with the collision so as not to grow in amplitude.  The
background radiation is invisible to the eye in solution A, and takes
the form of small-amplitude, non-localized, counter-propagating waves
in solution B when viewed as a movie at closer range than shown in
Figure~\ref{fig:cap:525}.  A Floquet stability analysis summarized in
Figures~\ref{fig:cap:705} and~\ref{fig:cap:525} shows that solution A
is linearly stable to harmonic perturbations while solution B is
unstable.

\section{Conclusion}

We have developed an efficient algorithm for computing the stability
spectra of standing water waves in finite and infinite depth, with and
without surface tension. The method involves two levels of truncation.
First, a truncated version of the monodromy operator is computed in
Fourier space by considering perturbations of the initial condition up
to a given wave number. This decouples the size of the truncated
operator from the size of the mesh, allowing all the matrix entries to
be computed accurately. Otherwise our method is similar to that of
Mercer and Roberts \cite{mercer:92}, who also solve the linearized
Euler equations about the time-periodic solution with multiple initial
conditions in order to study stability.  We then introduce a ``mean
wave number'' to order the eigenvalues of the truncated operator. This
ordering has the effect of moving accurately computed eigenvalues to
the front, as seen in plots of the residual error (\ref{eq:err})
obtained by substituting an eigenpair of the truncated operator into
the full state transition matrix after zero-padding the eigenvector.
We retain a specified number of Floquet multipliers (ordered in this
  way) and discard the rest, which is the second truncation step.
Finally, we define a linear assignment problem to match eigenvalues of
nearby standing waves in order to track individual eigenvalues from
the zero-amplitude wave to large-amplitude waves via homotopy. This
reveals, for example, which modes are involved when bubbles of
instability nucleate.

In addition to studying the stability of classical standing waves in
deep water, we explore the stability of new or recently discovered
families of time-periodic water waves that involve counter-propagating
solitary wave collisions, e.g.~of gravity waves in shallow water and
capillary-gravity waves in deep water. The examples in
Figures~\ref{fig:stable:005}AB, \ref{fig:concus} and~\ref{fig:cap:705}
show that large-amplitude water waves of various types can be stable
to harmonic perturbations both forward and backward in time. Computing
these waves in the first place is a significant challenge. One novelty
of the gravity-capillary waves computed in Figures~\ref{fig:cap:705}
and~\ref{fig:cap:525} is that we searched directly for large-amplitude
time-periodic solutions without using numerical continuation to get
there. This suggests an abundance of different types of time-periodic
water waves and also speaks to the robustness of the shooting method
to solve two-point boundary value problems with fairly inaccurate
starting guesses. False positives are avoided by resolving the
solutions with spectral accuracy and formulating the problem as an
overdetermined minimization problem.  In previous work
\cite{collision}, the author has faced similar challenges in computing
relative-periodic elastic collisions of co-propagating solitary water
waves that resemble cnoidal solutions of KdV.

Many of the numerical results of Section~\ref{sec:results} call for
future work. For example, what is the physical mechanism responsible
for the bubbles of instability we discovered with
crest acceleration far below the stability threshold $A_c=0.889$?
In a follow-up paper \cite{water:LT}, we will show that
perturbing the solution of the nonlinear equations near one of these
bubbles of instability leads to a cyclic pattern in which the
perturbation grows for thousands of cycles and then, surprisingly,
decays for thousands of cycles to return close to the standing wave in
a nearly time-reversed fashion to its initial excursion.  This
sequence repeats with some variation in how close the wave returns to
the standing wave before latching onto another growing mode. Bryant
and Stiassnie \cite{bryant:stiassnie:94} observed similar recurrent
behavior from sideband instabilities on a domain containing 9 replicas
of the standing wave, so this type of recurrent pattern may be a
common mechanism associated with a single pair of unstable Floquet
multipliers $\lambda$ and $\bar\lambda$.

Another natural question is the extent to which linear stability
predicts the long-time dynamics of perturbations of the nonlinear
equations. It has been observed in the literature
\cite{chan:street:70, maxworthy:76, su:mirie, cooker:97, 
  craig:guyenne:06, milewski:11} that although solitary water
wave collisions are inelastic, the residual radiation of such a
collision can be remarkably small. In \cite{water:LT}, we will show
that if two identical counter-propagating traveling water waves
(i.e.~Stokes waves) of a certain amplitude are combined at $t=0$ via
(\ref{eq:counter:prop:trav}), the solution remains close to solution B
of Figures~\ref{fig:bif:005}--\ref{fig:stable:005} for all 5000 cycles
we computed.  The first several collisions of the Stokes waves produce
radiation that becomes visible in the wave troughs after a few cycles;
however, it quickly saturates and does not grow beyond the size of the
background oscillations present in Solution B, which we found to be
linearly stable to harmonic perturbations. The amplitude of the
component traveling waves of the initial condition
(\ref{eq:counter:prop:trav}) was chosen to eliminate temporal drift
relative to solution B. This drift is caused by the secular
instability of one of the Jordan chains of length two associated with
$\lambda=1$. The secular instability of the other Jordan chain
was eliminated by even symmetry of the initial condition.

We were initially puzzled by this second Jordan chain. As explained in
Section~\ref{sec:jchains}, one of these chains corresponds to the
change in period when the amplitude of the standing wave is increased.
Equation (\ref{eq:jchain:eps}) shows that a perturbation in this
direction will cause the wave to drift forward or backward in time
relative to the underlying standing wave. We eventually realized that
the other Jordan chain corresponds to a perturbation direction in
which the wave drifts left or right in space as it evolves through
successive cycles. This suggests a bifurcation direction in which the
standing wave begins to travel as it oscillates.  The only reference
to such waves we have seen in the literature is a parenthetical
comment by Iooss, Plotnikov and Toland \cite{iooss05}: ``it is
possible to imagine more general solutions, for example,
`travelling-standing-wave' solutions, of the free boundary problem.''
Using the Jordan chains as a guide, we have successfully computed a
two-parameter family of traveling-standing waves that return to a
spatial translation of themselves at periodic time-intervals.  Pure
traveling waves and standing waves are special cases that occur at
certain values of the second bifurcation parameter.  These will be
reported on elsewhere \cite{waterTS}.

We have focused on harmonic stability in this paper, which is
natural for standing waves evolving in a rigid container. In future
work \cite{waterTS}, we will generalize the method to investigate
stability with respect to subharmonic perturbations in the more
general context of traveling-standing waves. Mercer and Roberts
\cite{mercer:92} considered subharmonic perturbations by replicating
the underlying standing wave and looking at harmonic stability over
a larger spatial period. Instead, we will consider quasi-periodic
perturbations of the linearized equations, similar to what has been
done recently for traveling waves \cite{benjamin67, benjamin67b,
  lh:stab:78a, lh:stab:78b, crawford:81, mackay:86, ioualalen:93,
  nicholls09, oliveras11, wilton:16}. There are many technical
challenges for the traveling-standing problem that do not arise in
the pure traveling case.

Subharmonic stability involves solving the linearized Euler equations
with quasi-periodic perturbations. A natural next step is to seek
fully nonlinear quasi-periodic solutions of the free-surface Euler
equations (\ref{eq:ww}). The relative-periodic solutions computed in
\cite{collision} and the traveling-standing waves computed in
\cite{waterTS} are examples of waves with two quasi-periods. A recent
paper of Berti and Montalto \cite{berti:montalto} proves existence of
gravity-capillary water waves with more than two quasi-periods using
Nash-Moser theory. We hope to extend the overdetermined shooting
method to compute such waves with or without surface tension and study
their properties. Spatially quasi-periodic traveling waves and other
generalizations also appear to be within reach with the tools we
are developing.

\appendix

\section{Boundary integral formulation}
\label{sec:BI}

In this section, we briefly describe how we compute the
Dirichlet-Neumann operator in (\ref{eq:DNO:def}).  Further details
(with derivations) may be found in Wilkening and Yu \cite{water2}. The
signs in (\ref{eq:KG}) below correct a typo in \cite{water2}, where
$K_1+K_2$ and $G_1+G_2$ were written in several equations without
changing the sign in the formulas for $K_2$ and $G_2$. The derivations
in \cite{water2} are otherwise correct.


To compute $\mc G\varphi$, we solve the integral equation
\begin{equation}\label{eq:fred}
  \frac{1}{2}\mu(\alpha) + \frac{1}{2\pi}\int_0^{2\pi}
  K(\alpha,\beta)\mu(\beta)\,d\beta =
  \varphi(\xi(\alpha))
\end{equation}
for the dipole density $\mu(\alpha)$, and then compute
\begin{equation}\label{eq:G}
    \mc G\varphi(\xi(\alpha)) = 
    \frac{1}{|\xi'(\alpha)|}
  \left[\frac{1}{2}H\gamma(\alpha) + \frac{1}{2\pi}\int_0^{2\pi}
    G(\alpha,\beta)\gamma(\beta)\,d\beta\right],
\end{equation}
where $\gamma(\alpha)=\mu'(\alpha)$ is the vortex sheet strength, $H$
is the Hilbert transform (with symbol $\hat H_k=-i\opn{sgn}(k)$), and
\begin{equation}
  x=\xi(\alpha)
\end{equation}
is a monotonic parametrization of the interval $[0,2\pi]$. Normally
$\xi(\alpha)=\alpha$, but other choices are useful for refining the
mesh in regions of high curvature \cite{water2}.  In these
formulas,
\begin{equation}\label{eq:KG}
  \begin{aligned}
    K(\alpha,\beta) &= K_1(\alpha,\beta) - K_2(\alpha,\beta), \\
    G(\alpha,\beta) &= G_1(\alpha,\beta) - G_2(\alpha,\beta),
  \end{aligned}
\end{equation}
where
\begin{alignat*}{2}
 K_1 &= \im\left\{\frac{\zeta'(\beta)}{2}
  \cot\frac{\zeta(\alpha) - \zeta(\beta)}{2}
  - \frac{1}{2}\cot\frac{\alpha-\beta}{2}\right\}, & \quad
 K_2 &= \im\left\{\frac{\bar\zeta'(\beta)}{2}\cot
  \frac{\zeta(\alpha)-\bar\zeta(\beta)}{2}\right\}, \\
 G_1 &= \re\left\{\frac{\zeta'(\alpha)}{2}
  \cot\frac{\zeta(\alpha) - \zeta(\beta)}{2}
  - \frac{1}{2}\cot\frac{\alpha-\beta}{2}\right\}, & \quad
 G_2 &= \re\left\{\frac{\zeta'(\alpha)}{2}\cot
  \frac{\zeta(\alpha)-\bar\zeta(\beta)}{2}\right\}.
\end{alignat*}
Here $\mathbb{R}^2$ has been identified with $\mathbb{C}$ and
the free surface is parametrized by
\begin{equation}\label{eq:zeta}
  \zeta(\alpha) = \xi(\alpha) + i\eta(\xi(\alpha)).
\end{equation}
When the fluid depth is infinite, $K_2$ and $G_2$ are dropped;
otherwise, the bottom boundary is assumed to be at $y=0$ so that
$\bar\zeta$ is the mirror image of the free surface.  
When $\beta=\alpha$, we set
\begin{align}
  K_1(\alpha,\alpha) &=
  -\im\{\zeta''(\alpha)/[2\zeta'(\alpha)]\}, \\
  G_1(\alpha,\alpha) &= \phantom{-}
  \re\{\zeta''(\alpha)/[2\zeta'(\alpha)]\},
\end{align}
which makes them continuous functions.

The integrals in (\ref{eq:fred}) and (\ref{eq:G}) are approximated
with spectral accuracy via the trapezoidal rule with uniformly spaced
collocation points $\alpha_j=2\pi j/M$, $0\le j<M$.  The matrix
entries $K_{ij} = K(\alpha_i,\alpha_j)$ and $G_{ij} =
G(\alpha_i,\alpha_j)$ are computed in parallel on a GPU, which
involves $O(M)$ communication costs for $O(M^2)$ work.  This makes
evaluation of $\phi_x$ and $\phi_y$ in (\ref{eq:ww}) and
(\ref{eq:uv:from:G}) comparable in speed to a conformal mapping
approach \cite{dyachenko:96b, nachbin:04, milewski:11}, which
involves only differentiation and the Hilbert transform on the
right-hand side of the evolution equations. Conformal mapping methods
are not suitable for computing extreme standing waves as the mesh
points spread out in sharply crested regions where mesh refinement is
needed.  Our formulation assumes the wave profile remains
single-valued; however, it is straightforward to generalize to
overturning waves \cite{lh76, baker:82, krasny:86, mercer:92, hls94,
  mercer:94, smith:roberts:99, ceniceros:99, HLS01, baker10, hmodel}.

\section{Computation of Traveling Water Waves}
\label{sec:trav}

To search for the new type of time-periodic gravity-capillary water
waves presented in Section~\ref{sec:surf}, we used initial guesses for
the shooting method consisting of two traveling waves moving in
opposite directions, superposed linearly at $t=0$, as in
(\ref{eq:counter:prop:trav}). Our trust-region shooting method is
easily adapted to compute such traveling waves.  It is not as
efficient as replacing $\eta_t$ and $\varphi_t$ in (\ref{eq:ww}) by
$-c\eta_x$ and $-c\varphi_x$ and solving the equations directly by a
variant of Newton's method \cite{rienecker,
  schwartz:fenton,milewski:11,wilton:16}, but it is
quite robust and does not require writing a separate code.

To find traveling waves, we minimize the objective function
\begin{equation}\label{eq:f:trav}
  f = \frac12r^Tr, \quad
  r_{2j} = \frac{\eta\left(x_j,\frac{T}{M}\right) -
    \eta\left(x_{j-1},0\right)}{\sqrt{2M}}, \qquad
  r_{2j+1} = \frac{\varphi\left(x_j,\frac{T}{M}\right) -
    \varphi\left(x_{j-1},0\right)}{\sqrt{2M}},
\end{equation}
where $M$ is the number of gridpoints, $j$ runs from 0 to $M-1$,
$x_j=2\pi j/M$, and $x_{-1}=x_{M-1}$.  We assume $\eta(x,0)$ is even
and $\varphi(x,0)$ is odd, replacing the initial condition
(\ref{eq:init:stand}) with
\begin{equation}\label{eq:init:trav}
  \hat\eta_k = c_{2|k|-1}, \qquad
  \hat\varphi_k = \pm ic_{2|k|}, \qquad (k=\pm1,\pm2,\dots \;;\; |k|\le n/2),
\end{equation}
where $n$ is an even integer.  As before $c_1,\dots,c_n$ are real and
all other Fourier modes are zero, except for $\hat\eta_0$ in the
finite depth case.  In the formula for $\hat\varphi_k$, the minus sign
is taken if $k<0$ so that $\hat\varphi_{-k}
=\overline{\hat\varphi_k}$.  We again define $T=c_0$ so that
$c\in\mathbb{R}^{n+1}$, and we can add an extra equation of the
form (\ref{eq:eta:constr}) with $m=2M$ to impose that $\eta(a,0)$
have a given value at $a=0$ or $a=\pi$.

Note that $f$ measures the difference between the solution at time
$T/M$ and a spatial shift of the initial condition by one grid point,
where $M$ is the number of grid points --- we assume
$\xi(\alpha)=\alpha$ in (\ref{eq:zeta}) when computing traveling
waves.  This objective function will be zero if $(\eta,\varphi)$ is a
traveling wave and $T$ is the time required to travel from $x=0$ to
$x=2\pi$.  Because the waves only travel to the right by one grid
point, a small number of time-steps (usually one or two) are typically
required to evolve the solution; thus, the method is fast.  A more
conventional approach for computing traveling waves is to substitute
$\eta(x-ct)$, $\varphi(x-ct)$ into (\ref{eq:ww}) and solve the
resulting stationary problem (or an equivalent integral equation) by
Newton's method \cite{rienecker, chen80a, 
  schwartz:fenton, chandler:93, milewski:11, wilton:16}.

\section{Computation of the Jacobian and the state transition matrix}
\label{sec:J}

To compute time-periodic solutions of (\ref{eq:ww}), we minimize
$f=\frac{1}{2}r^Tr$ in (\ref{eq:f:phi}) using a variant \cite{water2}
of the Levenberg-Marquardt method for nonlinear least squares problems
\cite{nocedal}.  This approach requires computing the Jacobian,
$J=\nabla_c r$.  As shown below, this can be done efficiently by
solving the variational equation (\ref{eq:variational}) with multiple
right-hand sides.  We compute the Fourier representation $\hat E_T$ of
the state transition matrix in (\ref{eq:ET:def}) using the same
technique.

As above, we suppress $x$-dependence in the notation when convenient.
Let $q(t) = (\eta(t);\varphi(t))$ represent the solution of
(\ref{eq:ww}) and $\dot q(t) = (\dot\eta(t); \dot\varphi(t))$
represent a derivative with respect to the initial condition (not
  time).  In more detail, we define
\begin{equation}\label{eq:eps:def}
  \dot q(t) = \der{}{\veps}\bigg\vert_{\veps=0} q(t;\veps), \qquad
  \begin{aligned}
    q_t(t;\veps) &= F(q(t;\veps)), \\
    q(0;\veps) &= q_0 + \veps \dot q_0,
  \end{aligned}
\end{equation}
where $F(q)$ denotes the right-hand side of (\ref{eq:ww}).  Each
column of the Jacobian is computed by solving the variational
equation (\ref{eq:variational}) for $\dot q$ alongside (\ref{eq:ww})
for $q$:
\begin{equation}\label{eq:q:qdot}
  \der{}{t}
  \begin{pmatrix} q \\ \dot q \end{pmatrix} =
  \begin{pmatrix} F(q) \\ DF(q)\dot q \end{pmatrix}, \quad
  \begin{aligned}
    q(0) &= q_0 = (\eta_0,\varphi_0), \\
    \dot q(0) &= \dot q_0 = \partial q_0/\partial c_k.
  \end{aligned}
\end{equation}
From the formula $r$ in (\ref{eq:f:phi}), we have
\begin{equation}\label{eq:Jjk}
  J_{jk} = \der{r_j}{c_k} = \begin{cases}
    (1/4)\varphi_t(\alpha_j,T/4)/\sqrt{M}, & k=0, \\
    \dot\varphi(\alpha_j,T/4)/\sqrt{M}, & k\ge1.
  \end{cases}
\end{equation}
In practice, $\dot q$ is replaced by the matrix $\dot Q=[\dot
  q_{(k=1)},\dots,\dot q_{(k=n)}]$ to compute all the columns of $J$
(besides $k=0$) at once.  This allows re-use of the matrices $K$ and
$G$ in the Dirichlet-Neumann operator across all the columns of the
Jacobian, and streamlines the linear algebra to run at level 3 BLAS
speed.

To compute $\hat E_T$, the initial condition (\ref{eq:q:qdot}) for
column $k$ of $J$ is replaced by (\ref{eq:q0:dot}) for column $4k+k'$
of $\hat E_T$; the linearized solutions are evolved from $t=0$ to
$t=T$ (instead of $T/4$); and, instead of (\ref{eq:Jjk}), the real and
imaginary parts of $\dot\eta_j^\wedge(T)$, $\dot\varphi_j^\wedge(T)$
are extracted to obtain rows $4j-3,\dots,4j$ of $\hat E_T$. We evolve
all the columns (or large batches of columns) in parallel, which
dramatically decreases the time required to compute all the entries
of $\hat E_T$.

\section{A matching problem for plotting eigenvalues smoothly}
\label{sec:matching}

The optional ``step (5)'' of Algorithm~\ref{alg:fmult} involves
matching the eigenvalues at adjacent values of $A_c$ to track
individual eigenvalues via homotopy from the zero-amplitude state to
large-amplitude standing waves.  This turns out to be surprisingly
challenging. Our solution is to formulate the problem of extracting
the smoothest possible ``eigenvalue curves'' as a sequence of linear
assignment problems \cite{munkres}.

The data for this section is the output of steps (1)--(4) of
Algorithm~\ref{alg:fmult}, which we implemented in C++, for the $380$
standing waves listed in the first 3 rows of
Table~\ref{tab:params}. As noted in the table, we set $n=600$ or
$n=900$ for the leading submatrix $J_{1:n,1:n}$ of $\hat E_T$ in step
(2) of the algorithm. For each of these standing waves, a file is
created with the leading $n^*$ computed eigenvalues $\lambda_j =
|\lambda_j|e^{i\sigma_j}$, sorted by mean wave number. The file
contains $|\lambda_j|$, $\sigma_j$, $\la k\ra_j$, and the parity $p_j$
for each eigenvalue. The numbers $p_j$ are set to 1 for odd
eigenfunctions and 0 for even ones. We use $n^*=360$ in steps
(1)--(4), match them in step (5) as explained below, and then discard
down to $300$, the reported value of $n^*$ in
Table~\ref{tab:params}. This allows us to track the same set of 300
eigenvalues via homotopy from $A_c=0$ to $A_c=0.9539$ even though they
may not all remain among the 300 eigenvalues of smallest mean wave
number. Figure~\ref{fig:kk:00} of Section~\ref{sec:deep} shows that the mean wave numbers at the
boundary between the retained and discarded eigenvalues vary smoothly and increase slowly, in lockstep,
over most of the range $0\le A_c\le 0.9539$. Thus, the same
result would have been obtained over most of this range using
$n^*=300$ for steps (1)--(4).  However, there are two instability
bubbles at the far right of the plot in Figure~\ref{fig:kk:00} that
lead to spikes with $\la k\ra$ exceeding 80. Here other eigenvalues
exist with smaller mean wave number, but we discard them as they are
not connected via homotopy to the curves of smallest mean wave number
when $A_c$ is small.

\begin{figure}
\begin{center}
\includegraphics[width=.9\linewidth]{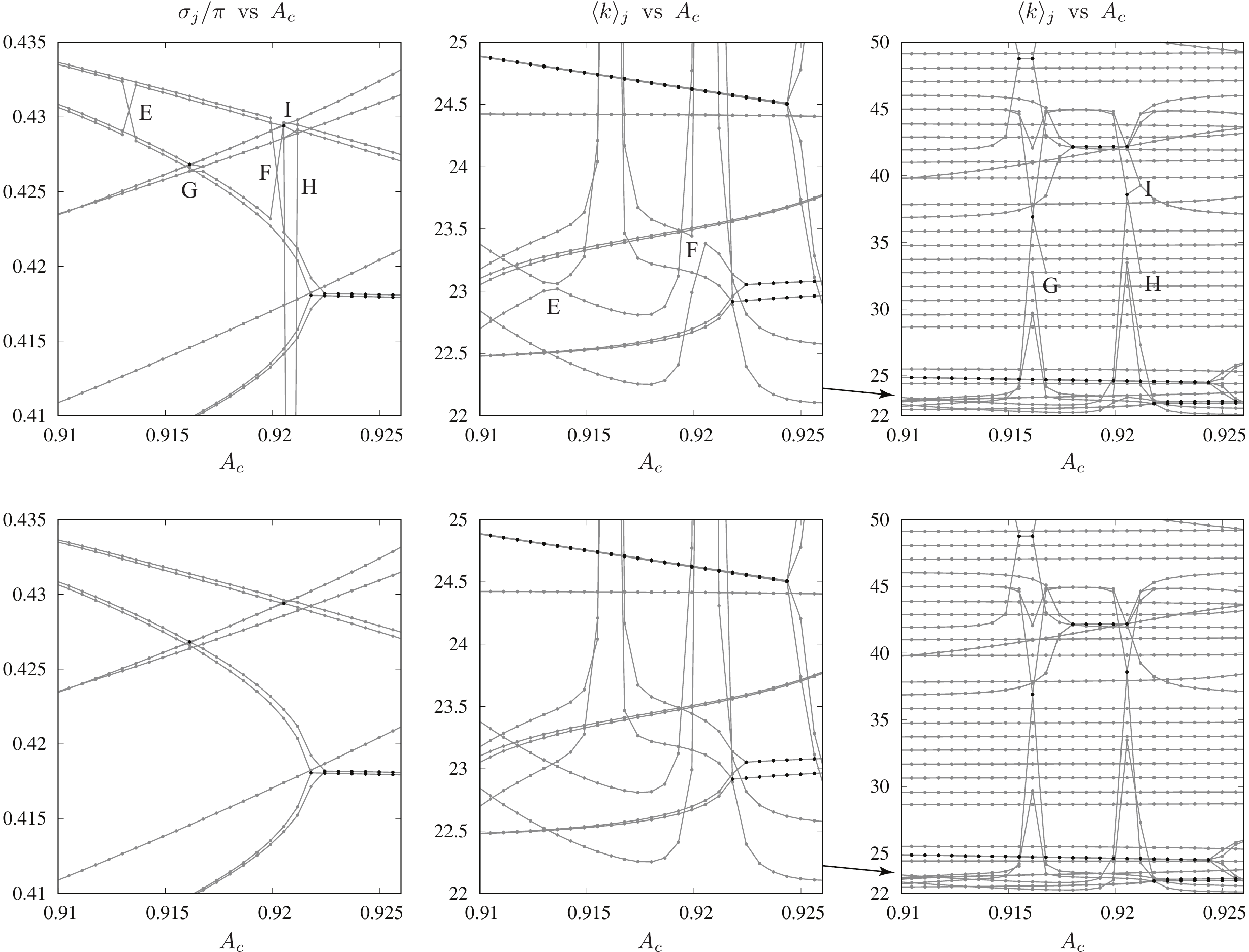}
\end{center}
\caption{\label{fig:match} If the cost matrix (\ref{eq:cost}) is
  modified by dropping the 1/2 powers and replacing the factor of 10
  by a factor of 2, errors arise in the matching algorithm (top row)
  that disappear when the cost function in (\ref{eq:cost}) is used
  (bottom row). Labels E, F, G, H, I show where the cost was reduced
  in $\la k\ra_j$ (top center and top right) at the expense of
  creating a criss-cross in $\sigma_j$ (top left).  }
\end{figure}

There are several matlab codes for the linear assignment problem
available in the public domain. We found the implementation described
in \cite{crouse:assign2D} to work well, which is based on the
algorithm of Jonker and Volgenant \cite{jonker:87}.  We combine the
data from the stability calculations into one large file and load it
into matlab to generate $n^*\times \ell$ matrices \verb+mag+,
\verb+arg+, \verb+kk+ and \verb+parity+, where $n^*=360$ and
$\ell=380$.  The columns of these matrices correspond to standing
waves with different values of crest acceleration $A_c$, sorted
smallest to largest. In particular, the first column of each of these
matrices contains the data for the $A_c=0$ wave.
Our goal is to generate an $n^*\times \ell$ matrix \verb|perm|,
defined so that \verb+i1=perm(i,s)+ means entry \verb|(i1,s)| of the
original matrices \verb+mag+, \verb+arg+, \verb+kk+ and \verb+parity+
should be moved to position \verb|(i,s)|. After re-ordering each
column, holding \verb|i| fixed as \verb+s+ runs from 1 to $\ell$ will
track a single eigenvalue through the family of standing waves. The
first column of \verb+perm+ is set to $[1:n^*]'$ so that the
eigenvalues at $A_c=0$, which we know analytically, remain ordered by
mean wave number.  Now suppose columns 1 through $s$ of \verb+perm+
have been computed and the entries of these columns in \verb+mag+,
\verb+arg+, etc.~have been permuted to their correct positions. The
linear assignment problem we propose is to find a permutation
$P=$\verb|perm(:,s+1)| of the integers $[1:n^*]$ to minimize the cost
function
\begin{equation}
  \sum_{i=1}^{n^*} C_{i,P(i)}^\e{s+1}.
\end{equation}
Through trial and error, we find the following cost matrix to be
effective:
\begin{equation}\label{eq:cost}
  \begin{split}
    C_{ij}^\e {s+1} &=
    10\left|\sigma_i^\e s + m_i^\e s\Big(A_c^\e{s+1}-A_c^\e s\Big) -
    \sigma_j^\e{s+1}\right|^{1/2} \\
    & \qquad\quad + \left|\big|\lambda_i^\e s\big| -
    \big|\lambda_j^\e {s+1}\big|\right|^{1/2}
    + \left|\langle k\rangle^\e s_i - \langle k\rangle^\e{s+1}_j\right|^{1/2}
    + 100\left|p_i^\e s - p_j^\e{s+1}\right|.
  \end{split}
\end{equation}
The slope $m_i^\e s$ used for linear extrapolation in (\ref{eq:cost})
is set to 0 if $\big|\lambda_i^\e s\big|$,
$\big|\lambda_i^\e{s-1}\big|$ or $\big|\lambda_j^\e{s+1}\big|$ differs
from 1 by more than $10^{-6}$, or if $s=1$.  Otherwise we define
\begin{equation}
  m_i^\e s = f\left(\frac{\sigma_i^\e s - \sigma_i^\e{s-1}}{
      A_c^\e s - A_c^\e{s-1}}\right),
  \quad
  f(x) = \begin{cases} 5 & x>5, \\ -5 & x < -5, \\ x & \text{o.w.}
    \end{cases}
\end{equation}
The $1/2$ powers in (\ref{eq:cost}) favor matching most eigenvalues
accurately and allowing a few to change significantly. This is
appropriate as eigenvalue collisions associated with bubbles of
instability cause the $|\lambda_j|$, $\sigma_j$ and $\la k\ra_j$
involved in the collision to change rapidly with $A_c$ while the other
eigenvalues change slowly.  The top row of Figure~\ref{fig:match}
shows the matching errors that appear if the 1/2 powers are replaced
by 1, and if the factor of 10 in front of the first term in
(\ref{eq:cost}) is replaced by 2. The bottom row shows the correct
results obtained via the cost matrix in (\ref{eq:cost}).  Most of the
curves turn out the same for the two choices of cost matrices, but
there are a handful of clear mistakes in the top row, labeled E, F, G,
H, I. The errors at E and F occur because short-circuiting the $\la
k\ra_j$ curves at E and F (rather than having them cross) in the top
middle panel reduces the cost in the modified version of
(\ref{eq:cost}) more than the cost increase introduced by crossing the
$\sigma_j$ curves at E and F in the top left panel. The errors at G,
H, I are due to large jumps in $\la k\ra_j$ associated with bubbles of
instability being too expensive in the modified cost function.
Introducing the 1/2 powers prevents the algorithm from breaking up
these large jumps into smaller jumps in $\la k\ra_j$ at the expense of
introducing erroneous crossings in the $\sigma_j$ curves.  With the cost
function (\ref{eq:cost}), we did not find a single connection error
anywhere in the data.

In our initial experiments trying different formulas for
$C_{ij}^\e{s+1}$ in (\ref{eq:cost}), there were often a few
criss-crossed curves such as in the top panels of
Figure~\ref{fig:match}. These are easily fixed by hand (outside of
  matlab). We write the matrix \verb|perm| to a file and create a
second text file containing 3 integers per line, \verb|s,i1,i2|. We
wrote a short perl script to read the permutation data into memory and
swap rows \verb|i1| and \verb|i2| in columns \verb|s| through $\ell$
of \verb|perm|. This is done repeatedly, once for each line in the
second file. Finally, the perl script reads the original eigenvalue
file (the file loaded by matlab to create \verb|mag|, \verb|arg|,
  \verb|kk| and \verb|parity|), re-orders the data according to the
corrected values in \verb|perm|, and writes a file with all the
eigenvalues ordered properly.  Finding the triples \verb|s,i1,i2|
required to correct errors such as at E, F, G, H or I in the top row
of Figure~\ref{fig:match} boils down to figuring out the indices
\verb|i1| and \verb|i2| associated with the two curves that cross incorrectly, and
the standing wave index \verb|s| where the crossing occurs.  This is
easily done by hand via a bisection algorithm. For each curve that
incorrectly crosses to another branch of eigenvalues, we use
gnuplot to plot ranges of indices in a different color, noting whether
the curve in question changes color or not. Repeatedly cutting the
number of curves plotted in half allows us to quickly identify the
index of the desired curve.  As mentioned above, the parameters in
(\ref{eq:cost}) yield correct matchings with no need for further
corrections to be performed by hand.

\bibliographystyle{plain}

\end{document}